\documentclass[a4paper, 11pt]{article}

\usepackage[utf8]{inputenc}
\usepackage{lineno}
\usepackage{bm}
\usepackage{fullpage}
\usepackage{graphicx}
\usepackage{psfrag}
\usepackage{amsmath}
\usepackage{amsfonts}
\usepackage{bbm}
\usepackage{bbold}
\usepackage{verbatim}
\usepackage[colorlinks=false]{hyperref}
\usepackage{enumerate}
\usepackage[shortlabels]{enumitem}
\usepackage{subfig}
\usepackage{algorithm}
\usepackage{algpseudocode}
\usepackage[noabbrev]{cleveref}
\usepackage{soul}
\usepackage{tabularx}
\usepackage{booktabs}
\usepackage{multicol}
\usepackage{multirow}
\usepackage{caption}
\usepackage{natbib}
\usepackage{xcolor}

\setcitestyle{authoryear,round}

\title{Dynamically adaptive networks for integrating optimal pressure management and self-cleaning controls}

\author{Bradley Jenks$^1$, Aly-Joy Ulusoy$^1$, Filippo Pecci$^{1, 2}$, Ivan Stoianov$^1$}
\date{
    $^1$InfraSense Labs, Department of Civil and Environmental Engineering, Imperial College London, London, SW7 2BB, UK\\
    $^2$Andlinger Center for Energy and the Environment, Princeton University, Princeton, New Jersey, 08544, United States
}

\begin{document}
\maketitle
\begin{abstract}
    This paper investigates the problem of integrating optimal pressure management and self-cleaning controls in dynamically adaptive water distribution networks. We review existing single-objective valve placement and control problems for minimizing average zone pressure (AZP) and maximizing self-cleaning capacity (SCC). Since AZP and SCC are conflicting objectives, we formulate a bi-objective design-for-control problem where locations and operational settings of pressure control and automatic flushing valves are jointly optimized. We approximate Pareto fronts using the weighted sum scalarization method, which uses a previously developed convex heuristic to solve the sequence of parametrized single-objective problems. The resulting Pareto fronts suggest that significant improvements in SCC can be achieved for minimal trade-offs in AZP performance. Moreover, we demonstrate that a hierarchical design strategy is capable of yielding good quality solutions to both objectives. This hierarchical design considers pressure control valves first placed for the primary AZP objective, followed by automatic flushing valves placed to augment SCC conditions. In addition, we investigate an adaptive control scheme for dynamically transitioning between AZP and SCC controls. We demonstrate these control challenges on case networks with both interconnected and branched topology.
\end{abstract}

\definecolor{royal_blue}{rgb}{0.25, 0.41, 0.88}
\definecolor{orange_yellow}{rgb}{1.0, 0.65, 0.0}


\section{Introduction}
\label{sec:intro}
The optimal design and control of water distribution networks (WDNs) is critical to maintain performance  under increasing environmental, regulatory, and financial constraints. Advanced control structures and methods aim to satisfy multiple operational criteria, many of which are conflicting objectives with regard to network hydraulics. In this paper, we investigate the bi-objective problem of integrating optimal pressure management and self-cleaning controls. These are important operational objectives linked to the reduction of leakage and discolouration risk in WDNs, respectively. However, high complexity and insufficient control capabilities have precluded the integration of pressure and water quality operational control schemes. Water companies are therefore seeking insights to inform design and control interventions for optimally integrating water quality as part of an overall control scheme.

Pressure management is a well-researched and implemented operational practice to lessen both the extent of background leakage and burst frequency in WDNs \citep{LAMBERT2010,SCHWALLER2015}. Typical pressure management schemes consist of pressure control valves (PCVs) installed at the inlets to hydraulically isolated areas of a WDN. This type of network configuration is called a pressure managed area (PMA), in which PCV actuators modulate downstream pressure and flow is continuously monitored. Note that a PMA coincides with the practice of district metered areas (DMAs), but with pressure control capabilities. Since leak flow is a function of hydraulic pressure, the implementation of PMAs is a cost-effective measure for reducing background leakage \citep{PILLER2014,SCHWALLER2015}. For example, PMAs have been successfully implemented by the UK water industry in the last 25 years, acting as a key strategy within leakage management programmes \citep{OFWAT2022}. The problem of optimal design and control of PMAs has therefore been the focus of numerous studies in recent literature. These have commonly been formulated as minimization problems with average zone pressure (AZP) used as a surrogate for pressure induced leak flow \citep{WRIGHT2015}. We present a review of such literature in \Cref{sec:review_obj_pressure}. In practice, however, pressure management schemes often employ a relatively simplistic fixed-outlet mode of pressure control. This, in conjunction with kept-shut boundary valves, permanently sectorizes the network. One of many potential consequences of this pressure management scheme is the deterioration of water quality, often observed as an increase in discolouration incidents \citep{MACHELL2014,WRIGHT2014,ARMAND2018}. 

Discolouration incidents are the leading cause of customer complaints in WDNs \citep{ARMAND2017}. Although a growing body of research has focused on the physical mechanisms leading to their occurrence \citep[e.g.][]{BOXALL2005,VANSUMMEREN2017,ARMAND2018,BRAGA2020}, little proactive maintenance is currently performed in operational networks. Existing methods for mitigating discolouration risk include manual flushing programmes, during which water is flushed from a network through hydrants. Such flushing activity is resource-intensive and often performed in response to customer complaints. In addition to flushing being a costly operational activity, increasingly strict water quality regulations are making compliance measures more challenging for water companies. This is especially the case in the UK water industry, where the Drinking Water Inspectorate (DWI) enforce financial penalties for failure to comply with water quality performance targets \citep{DWI2020}. Water companies are therefore seeking cost-effective and automatic control strategies to reduce the severity and frequency of discolouration risk. A control strategy of growing interest in both practice and academic literature is the concept of self-cleaning networks. This leverages the intrinsic connection between discolouration processes and hydraulic conditions \citep[see review in][]{ARMAND2018} by controlling peak diurnal flow velocities. More specifically, the strategy seeks to increase the distribution of velocities above a self-cleaning threshold so as to periodically re-suspend and flush accumulated material through the network \citep{VREEBURG2009,BLOKKER2012}. Recent studies have formulated this as an optimization problem, with decision variables comprising the placement of isolation valves \citep{ABRAHAM2018}, the operational settings of PCVs \citep{ABRAHAM2016}, and the placement and operational settings of PCVs and automatic flushing valves (AFVs) \citep{JENKS2023}. Note that AFVs are mechanical actuators which facilitate a controlled operational demand for flushing the network at hydrants or dedicated installations. These are increasingly deployed devices in operational networks and represent a key control actuator for improving SCC performance. We review the literature on optimizing SCC conditions in \Cref{sec:review_obj_self-cleaning}.

To facilitate the inclusion of multiple control objectives, \citet{WRIGHT2014,WRIGHT2015} introduced the operational framework of dynamically adaptive WDNs. This framework enables the dynamic re-configuration of network connectivity and hydraulic control to modulate flow and pressure. In this paper, we investigate the problem of integrating optimal pressure management and self-cleaning controls through the design and operation of dynamically adaptive networks. To the best of the authors’ knowledge, the AZP-SCC bi-objective problem has not been considered in previous literature.

Following a review of existing single-objective control strategies in \Cref{sec:review}, we present the two main contributions of this work. First, we propose a bi-objective framework for simultaneously considering AZP and SCC objectives. The bi-objective problem extends the work suggested in \citet{ABRAHAM2019} by jointly optimizing the placement and operation for both PCV and AFV actuators. This includes an investigation of a hierarchical design strategy, which first considers the placement of PCVs for the AZP single-objective problem, followed by AFVs for the SCC single-objective problem. We use the weighted sum scalarization method to formulate a sequence of parameterized single-objective optimization problems for approximating the AZP-SCC Pareto front. The problem formulation and solution methodology are presented in \Cref{sec:bi_objective_formulation}. Numerical experiments using two case networks demonstrate single-objective AZP and SCC results in \Cref{sec:num_experiment_single_dfc} and trade-offs between objectives through approximated Pareto fronts in \Cref{sec:num_experiment_bi_dfc}. These results aim to inform decision making for the design of PCV and AFV actuators that optimally integrate AZP and SCC objectives. In the second contribution, we investigate an adaptive control problem to facilitate the dynamic transition between AZP and SCC operational modes. The adaptive control problem considers only the operational settings of PCV and AFV actuators, assuming their placement is selected on the basis of the bi-objective problem results. Numerical experiments in \Cref{sec:num_experiment_adaptive_control} demonstrate the challenges in handling large nodal pressure variations that arise during the transition between AZP and SCC controls. We conclude with suggestions for future work to constrain maximum pressure and diurnal pressure variation in the AZP-SCC adaptive control problem.

\section{Preliminaries and literature review}
\label{sec:review}
In this section, we review existing control strategies for optimal pressure management and self-cleaning operations. The reviewed approaches focus on computing the optimal locations and operational settings of pressure control and automatic flushing valve actuators. First, \Cref{sec:review_hydraulics} reviews the hydraulic variables and conservation equations used to define the feasible solution space. We then present the motivation and mathematical formulations for pressure management and self-cleaning objective functions in \Cref{sec:review_obj_pressure} and \Cref{sec:review_obj_self-cleaning}, respectively. Following this, \Cref{sec:review_minlp} describes the general single-objective valve placement and control problem formulation. Thereafter, a brief review of previously implemented solution methods for single-objective problems is included in \Cref{sec:review_solution}. On the basis of this review, we identify opportunities and challenges for the integration of pressure management and self-cleaning control strategies in \Cref{sec:review_gaps}.

\subsection{Hydraulic variables and constraints}
\label{sec:review_hydraulics}
We model a water distribution network (WDN) as a directed graph $G(V, E)$ with $|V|$ vertices comprising $n_n$ demand nodes and $n_0$ source nodes (e.g. reservoirs), and $|E|$ edges comprising $n_p$ links. Network connectivity is represented through link-node incidence matrices $A_{12} \in \mathbb{R}^{n_p \times n_n}$ and $A_{10} \in \mathbb{R}^{n_p \times n_0}$ for demand and source nodes, respectively. The following convention is used to define $A_{12}$ \citep{TODINI1988}: 
\begin{linenomath}
\begin{equation} \label{eq:incidence_matrix_A12}
    \begin{alignedat}{3}
        &(A_{12})_{ji} =
        \begin{cases} 
            1 &\text{if link $j$ enters node $i$} \\
            0 & \text{if link $j$ is not connected to node $i$} \\
            -1 &\text{if link $j$ leaves node $i$.}
        \end{cases}
    \end{alignedat}
\end{equation}
\end{linenomath}
The same convention is used to define incidence matrix $A_{10}$ for source nodes. Since link directions are defined \textit{a priori}, any flows computed in the opposite direction are simply assigned a negative value.

We consider a steady-state hydraulic analysis over $n_t$ discrete time steps in an extended period simulation (EPS). For each time step ${t \in \{1, \dots, n_t\}}$, known hydraulic conditions are given by vectors of nodal demands $d_t \in \mathbb{R}^{n_n}$ and source hydraulic heads $h^0_t \in \mathbb{R}^{n_0}$. Moreover, we include vectors $\eta_t \in \mathbb{R}^{n_p}$ and $\alpha_t \in \mathbb{R}^{n_n}$ to model local losses introduced across pressure control valves (PCVs) and operational demands at automatic flushing valves (AFVs), respectively. Unique vectors of hydraulic states $q_t \in \mathbb{R}^{n_p}$ and $h_t \in \mathbb{R}^{n_n}$ are computed by solving the following energy \eqref{eq:hyd_energy_nonlin}-\eqref{eq:hyd_energy_lin} and mass \eqref{eq:hyd_mass} conservation equations:
\begin{linenomath}
\begin{subequations}
    \begin{align}
        &\theta_t - \phi(q_t) = 0, \label{eq:hyd_energy_nonlin} \\
        &A_{12}h_t + A_{10}h^0_t + \theta_t + \eta_t = 0,
        \label{eq:hyd_energy_lin} \\
        &A_{12}^T q_t - d_t - \alpha_t = 0, \label{eq:hyd_mass}
    \end{align}
\end{subequations}
\end{linenomath}
where the vector of auxiliary variables $\theta_t \in \mathbb{R}^{n_p}$ is introduced to isolate the nonlinear head loss model $\phi(q_t)$. The vector $\phi(q_t) = [\phi_1(q_{1,t}) \dots \phi_{n_p}(q_{n_{p},t})]^T$ is defined in general form as
\begin{linenomath}
\begin{equation}\label{eq:head_loss_model}
    \phi_j(q_{jt}) = r_j|q_{jt}|^{n_j-1}q_{jt}, \quad \forall j \in \{1,\ldots,n_p\}.
\end{equation}
\end{linenomath}
The resistance coefficient $r_j$ and exponent $n_j$ take different values depending on the link type (e.g. pipe or valve) and on the head loss model used. For valve links, $n_j = 2$ and 
\begin{linenomath}
\begin{equation} \label{eq:local_loss_model}
    r_j = \frac{8K_j}{g\pi^2D_j^4},
\end{equation}
\end{linenomath}
where $K_j$ and $D_j$ are the valve loss coefficient and diameter, respectively \citep{LAROCK1999}. For frictional head losses across pipe links, the Hazen-Williams (HW) or Darcy-Weisbach (DW) formulae are commonly applied. In the HW model, $n_j = 1.852$ and
\begin{linenomath}
\begin{equation} \label{eq:HW_model}
    r_j = \frac{10.67L_j}{C_j^{1.852}D_j^{4.871}},
\end{equation}
\end{linenomath}
where $C_j$ is the HW coefficient, a dimensionless number representing pipe frictional characteristics; $L_j$ is pipe length; and $D_j$ is pipe diameter \citep{LAROCK1999}. When the DW formula is used, there exists an implicit relationship between $q_{jt}$ and $r_j$, resulting in a continuously differentiable function for $\phi(q_{jt})$. The reader is referred to \citet{WALDRON2020} for details on formulating the DW model. While these models are ubiquitous in pipeline hydraulics, they both comprise nonsmooth terms (e.g. unbounded second order derivatives) and thereby present challenges for the use of state-of-the-art nonlinear optimization solvers. A smooth quadratic approximation (QA) for either the HW or DW model has been adopted in many WDN optimization studies \citep[e.g.][]{ECK2013,PECCI2019,ULUSOY2020}. These approximations follow the methodologies presented in \citet{ECK2015} and \citet{PECCI2017a}, where a pair of coefficients $(a_j,b_j)$ are generated for each link to form the QA head loss model,
\begin{linenomath}
\begin{equation} \label{eq:QA_model}
    \phi(q_{jt}) = q_{jt}(a_j|q_{jt}|+b_j).
\end{equation}
\end{linenomath}
Although the QA model is convenient for application in nonlinear optimization problems, it introduces errors in the hydraulic variables $q_t$ and $h_t$. In many cases, these errors are small in comparison to the degree of uncertainty inherent in hydraulic models and thus do not materially impact the feasibility of the optimal solutions \citep{GHADDAR2017,PECCI2017a, ULUSOY2022b}. However, for optimization problems which aim to increase flows, such as the self-cleaning capacity problem \citep{ABRAHAM2016,ABRAHAM2018,JENKS2023}, the large range of flow values over which the QA model is defined may result in errors \citep[Equation ~21]{PECCI2017a}. Such QA errors need to be quantified and assessed to ensure they are within an acceptable tolerance.

The placement of valve actuators is modelled as follows. Binary variables $z \in \{0,1\}^{n_p}$ model the placement of PCVs at network links. For each time step $t \in \{1, \dots, n_t\}$, $v^+_t \in \{0,1\}^{n_p}$ and $v^-_t \in \{0,1\}^{n_p}$ assign pressure control capabilities in the positive or negative flow direction, respectively. Control is limited to a single direction at each time step $t$ by the following constraint:
\begin{linenomath}
\begin{equation} \label{eq:valves_physical_a}
    v^+_{jt} + v^-_{jt} \leq z_j, \quad \forall  j \in \{1, \dots, n_p\},\,\; \forall t \in \{1, \dots, n_t\}
\end{equation}
\end{linenomath}
Moreover, binary variables $y \in \{0,1\}^{n_n}$ model the placement of AFVs at network nodes. Binary variables $z$ and $y$ are subject to the following physical constraints, which enforce the number of installed pressure control $n_v$ and flushing $n_f$ valves:
\begin{linenomath}
\begin{subequations}
\begin{align}
    &\sum_{j=1}^{n_p}z_j = n_v, \label{eq:valves_physical_b}\\
    &\sum_{i=1}^{n_n}y_i = n_f. \label{eq:valves_physical_c}
\end{align}
\end{subequations}
\end{linenomath}

We introduce constant vectors to define lower and upper bounds on continuous variables $q_t$, $h_t$, $\eta_t$, $\theta_t$, and $\alpha_t$, such that, for all $t \in \{1,\ldots,n_t\}$,
\begin{linenomath}
\begin{subequations}
\label{eq:hyd_bounds}
\begin{align}
    &q^L_t \leq q_t \leq q^U_t, \label{eq:hyd_bounds_a}\\
    &h^{\min}_t \leq h_t \leq h^{\max}_t, \label{eq:hyd_bounds_b}\\
    &\eta^L_t \leq \eta_t \leq \eta^U_t, \label{eq:hyd_bounds_c}\\
    &\theta^L_t \leq \theta_t \leq \theta^U_t, \label{eq:hyd_bounds_d}\\
    &0 \leq \alpha_t \leq \alpha^U_t, \label{eq:hyd_bounds_e}
\end{align}
\end{subequations}
\end{linenomath}
where $q^L_t$ and $q^U_t$ are derived from a given vector of maximum allowable velocities ${u^{\max} \in \mathbb{R}^{n_p}}$; $h^{\min}_t$ and $h^{\max}_t$ are set by minimum regulatory pressures and maximum known source heads, respectively; $\eta^L_t$ and $\eta^U_t$ are derived from corresponding $h^{\min}_t$ and $h^{\max}_t$ values at the upstream and downstream nodes of PCV links \citep[Equation ~3]{PECCI2019}; $\theta^L_t$ and $\theta^U_t$ are computed by \eqref{eq:hyd_energy_nonlin} with flow bounds $q^L_t$ and $q^U_t$; and $\alpha^U_t$ is defined with knowledge of local network conditions.

Variable bounds are modified for network links and nodes assigned as PCV and AFV actuators, respectively. We formulate \mbox{big-M} constraints to activate actuator settings $\eta_t$ and $\alpha_t$ and enforce energy conservation across PCV links. These constraints are written as follows, for all $t \in \{1,\ldots,n_t\}$:
\begin{linenomath}
\begin{subequations}
\label{eq:valve_bigM}
\begin{align}
    &\eta_t - \text{diag}(\eta^U_t)v^+_t \leq 0, \label{eq:valve_bigM_a}\\
    &-\eta_t + \text{diag}(\eta^L_t)v^-_t \leq 0, \label{eq:valve_bigM_b}\\
    &-q_t - \text{diag}(q^L_t)v^+_t \leq -q^L_t, \label{eq:valve_bigM_c}\\
    &q_t + \text{diag}(q^U_t)v^-_t \leq q^U_t, \label{eq:valve_bigM_d}\\
    &-\theta_t - \text{diag}(\theta^L_t)v^+_t \leq -\theta^L_t, \label{eq:valve_bigM_e}\\
    &\theta_t + \text{diag}(\theta^U_t)v^-_t \leq \theta^U_t, \label{eq:valve_bigM_f}\\
    &\alpha_t - \text{diag}(\alpha^U_t)y \leq 0. \label{eq:valve_bigM_g}
\end{align}
\end{subequations}
\end{linenomath}

In summary, we denote a feasible solution to the optimal control problem as the matrix $x = [x_1, \dots, x_{n_t}]^T$, where $x_t = [q_t, h_t, \eta_t, \theta_t, \alpha_t]$ are the continuous hydraulic variables for all $t \in \{1, \dots, n_t\}$. The vector $x_t$ satisfies hydraulic constraints \eqref{eq:hyd_energy_nonlin}-\eqref{eq:hyd_mass} and \eqref{eq:hyd_bounds_a}-\eqref{eq:valve_bigM_g}. Moreover, we denote optimal valve actuator placement through binary variables $z$, $v_t := [v^+_t, v^-_t]^T$, and $y$, which satisfy physical and economical constraints \eqref{eq:valves_physical_a} and \eqref{eq:valves_physical_b}-\eqref{eq:valves_physical_c}.

\subsection{Objective functions}
\label{sec:review_obj}
This section outlines the objective functions used in optimal pressure management and self-cleaning control problems. These represent surrogate measures for the actual operational conditions sought to be optimized.

\subsubsection{Pressure management}
\label{sec:review_obj_pressure}
The problem of optimal valve control for pressure management has been a subject of growing interest in the WDN optimization research community. This control problem seeks to reduce and maintain network pressures as close as possible to the minimum regulatory requirements in order to lessen the extent of background leakage and pipe circumferential stress \citep{SCHWALLER2015}. The network's average zone pressure (AZP) is commonly used as a surrogate for pressure induced leak flow, which is defined as follows \citep{WRIGHT2015}:
\begin{linenomath}
\begin{equation} \label{eq:azp_objective}
    \text{AZP}(x) := \frac{1}{n_t} \sum_{t=1}^{n_t} \sum_{i=1}^{n_n} w_i (h_{it} - \zeta_i),
\end{equation}
\end{linenomath}
where $h_{it}$ is the computed hydraulic head for node $i$ at time step $t$; $\zeta_i$ is the elevation of node $i$; and $w_i$ is a coefficient weighting node $i$ by the length of its connected links. That is,    
\begin{linenomath}
\begin{equation}
    \label{eq:azp_weights_a}
        w_i = \bar{L}^{-1} \sum_{k \in \mathcal{I}_i} \frac{L_i}{2},
\end{equation}
\end{linenomath}
where $\mathcal{I}_i$ is the set of links incident to node $i$; $\bar{L} = \sum_{j=1}^{n_p} L_j$; and $L \in \mathbb{R}^{n_p}$ is the vector of pipe lengths.

Minimizing AZP (or an analogous function) has been the objective of numerous optimal valve control problems. These have been formulated as single-objective problems \citep[e.g.][]{JOWITT1989,REIS1997,VAIRA1998,ULANICKI2000,ARAUJO2006,ECK2013,DAI2014,WRIGHT2014,WRIGHT2015,COVELLI2016,DEPAOLA2017,GHADDAR2017,PECCI2019,ULUSOY2020}, as part of multi-objective problems \citep[e.g.][]{NICOLINI2009,CREACO2015,PECCI2017b,NERANTZIS2020,ULUSOY2022a,ULUSOY2022b}, or included as joint optimization problems with other objectives \citep[e.g.][]{PECCI2021,PECCI2022a}. For multi-objective problems, we note that AZP has often been considered as the main operational objective. For example, \citet{ULUSOY2020} perform an \textit{a posteriori} analysis of network resiliency using solutions from a single-objective AZP problem. Moreover, \citet{NERANTZIS2022} propose an adaptive model predictive control (MPC) scheme where AZP is implemented as the main (or continuous) objective. 

In addition to reducing the mean operating pressure, recent studies have highlighted the influence of nodal pressure variation on the probability of pipe failure \citep{REZAEI2015,MARTINEZ2016,ARRIAGADA2021}. Accordingly, efforts have been made to introduce pressure variability (PV) as an objective within pressure management strategies. On the basis of PV indicators posed in \citet{MARTINEZ2016}, \citet{PECCI2017b} defined the PV objective function as the sum of squared differences between nodal heads $h_t$ and $h_{t-1}$, as follows:
\begin{linenomath}
\begin{equation} \label{eq:pv_objective}
    \text{PV}(x) := \sum_{i=1}^{n_n} (h_{i \,1}-h_{i \,n_t})^2 + \sum_{t=2}^{n_t} \sum_{i=1}^{n_n} (h_{it}-h_{i \,t-1})^2 ,
\end{equation}
\end{linenomath}

In demonstrating the impacts from optimal AZP valve control (i.e. flow modulation curves), \citet{WRIGHT2014} reported an increase in nodal pressure variability when compared with a fixed-outlet pressure management strategy. This suggested that AZP and PV are conflicting objectives and therefore cannot be simultaneously optimized. To this end, \citet{PECCI2017b} investigated trade-offs between AZP and PV through a bi-objective problem formulation. The authors applied different scalarization methods to generate the set of Pareto optima. In this work, we introduce PV as a critical metric to consider \textit{a posteriori} when implementing control schemes in operational WDNs. We discuss these challenges further in \Cref{sec:num_experiment_adaptive_control}.

\subsubsection{Self-cleaning capacity}
\label{sec:review_obj_self-cleaning}
Reducing discolouration risk is increasingly becoming a key priority in water quality management programmes. Since conventional discolouration risk reduction practices are infrequent and resource-intensive (e.g. manual flushing), there have been recent forays in the design and operation of self-cleaning networks \citep{BLOKKER2012}. Adapted from \citet{VREEBURG2009}, the self-cleaning capacity (SCC) is defined as the percentage of the network (by pipe length) routinely experiencing flow velocities above a threshold adequate to mobilize accumulated particles. Previous optimization studies have used SCC as a surrogate for discolouration risk, which is defined as follows \citep{ABRAHAM2016,ABRAHAM2018,JENKS2023}:
\begin{linenomath}
\begin{equation} \label{eq:scc_objective}
    \text{SCC}(x) := \frac{1}{n_t} \sum_{t=1}^{n_t} \sum_{j=1}^{n_p} w_j \kappa_j\bigg(\frac{q_{jt}}{A_j}\bigg),
\end{equation}
\end{linenomath}
where $q_{jt}$ is the computed flow conveyed by link $i$ at time step $t$; $A_j$ is the cross-sectional area of link $j$; and $\kappa_j(\cdot)$ is an indicator function which models the state of pipe velocities with reference to a minimum threshold. The indicator function described for link $j \in \{1, \dots, n_p\}$ as
\begin{linenomath}
\begin{equation}
\label{eq:scc_objective_indicator}
    \kappa_j(u) = 
    \begin{cases} 1 &\text{if $|u|$ $>$ $u^{\min}_j$} \\
    0 & \text{otherwise},
    \end{cases}
\end{equation}
\end{linenomath}
where $u^{\min}_j$ represents the minimum threshold velocity at link $j$ \citep[see][]{JENKS2023}. Moreover, a weighting $w_j$ normalizes the length of link $j$ to the entire network, $w_{j} = \frac{L_j}{\sum_{k=1}^{n_p}L_k}$.

Because \eqref{eq:scc_objective} is nonsmooth at $\pm u^{\min}$, \citet{ABRAHAM2016} proposed an approximation of $\kappa(\cdot)$ using a smooth sum of sigmoids (or logistic) function. The sum of sigmoids function is defined by positive and negative components, $\psi^+_j(u) := (1+e^{-\rho(u-u^{\min}_j)})^{-1}$ and $\psi^-_j(u) := (1+e^{-\rho(-u-u^{\min}_j)})^{-1}$, respectively, where parameter $\rho$ sets the sigmoid function curvature. A smooth approximation of \eqref{eq:scc_objective} is then defined as
\begin{linenomath}
\begin{equation}
\label{eq:scc_objective_sigmoid}
    \widetilde{\text{SCC}}(x) := \frac{1}{n_t} \sum_{t=1}^{n_t} \sum_{j=1}^{n_p}w_j \left(\psi^+_{j}\bigg(\frac{q_{jt}}{A_j}\bigg) + \psi^-_j\bigg(\frac{q_{jt}}{A_j}\bigg)\right),
\end{equation}
\end{linenomath}
to which gradient-based optimization methods can be applied.

Although water quality has been the objective of numerous studies in the WDN optimization literature \citep[see][]{MALA2017}, little research has focused on maximizing SCC to reduce discolouration risk. In \citet{ABRAHAM2016}, a single-objective optimization problem was formulated to maximize SCC through optimal control valve settings. Similarly, \citet{ABRAHAM2018} redistributed flow by optimizing isolation valve closures through a linear graph analysis tool. The developed algorithm was shown to outperform an empirical strategy where valve closures were based on engineering judgment. More recently, \citet{JENKS2023} formulated a novel design-for-control problem where locations and operational settings of PCV and AFV actuators were jointly optimized to maximize SCC. \Cref{sec:review_minlp} presents the general valve placement and control problem formulation.

\subsection{Valve placement and control problem formulation}
\label{sec:review_minlp}
The single-objective valve placement and control problem aims to minimize AZP (\Cref{sec:review_obj_pressure}) or maximize SCC (\Cref{sec:review_obj_self-cleaning}), subject to hydraulic modelling constraints and economical valve placement constraints (\Cref{sec:review_hydraulics}). The corresponding mixed integer nonlinear programming (MINLP) problem is written as 
\begin{linenomath}
\begin{subequations} \label{eq:minlp_problem}
    \begin{alignat}{3}
        & \underset{\substack{x, \, v, \, z, \, y}}{\text{minimize}}
        & \quad & f(x) \label{eq:minlp_problem_a} \\
        & \text{subject to}
        & & Ax - \phi(x) = 0 \label{eq:minlp_problem_b} \\
        & & & x \in \mathcal{X}(v,y), \; v \in \mathcal{V}(z)\label{eq:minlp_problem_c} \\
        & & & v \in \{0,1\}^{2n_pn_t}, \; z \in \{0,1\}^{n_p}, \; y \in \{0,1\}^{n_n},
    \end{alignat}
\end{subequations}
\end{linenomath}
which is described as follows. First, $f(x)$ in \eqref{eq:minlp_problem_a} denotes the objective function, set to either AZP \eqref{eq:azp_objective} or the additive inverse of SCC \eqref{eq:scc_objective_sigmoid} in this study. Matrix $A$ is defined so that the rows of \eqref{eq:minlp_problem_b} correspond to the vector of auxiliary variables $\theta_t \in \mathbb{R}^{n_p}$. Moreover, $\phi(\cdot)$ is a nonlinear function representing the frictional head loss model, as described in \Cref{sec:review_hydraulics}. The set $\mathcal{V}$ is defined by linear valve placement constraints \eqref{eq:valves_physical_a} and \eqref{eq:valves_physical_b}-\eqref{eq:valves_physical_c}. Given $v=(v^+,v^-) \in \mathcal{V}$, the polyhedral set $\mathcal{X}(v,y)$ is defined by linear hydraulic conservation constraints \eqref{eq:hyd_energy_lin}-\eqref{eq:hyd_mass}, big-M constraints for modelling valve actuator settings \eqref{eq:valve_bigM_a}-\eqref{eq:valve_bigM_g}, and lower and upper bounds on continuous variables \eqref{eq:hyd_bounds_a}-\eqref{eq:hyd_bounds_e}. The continuous decision variables are defined as: $q := (q_t)_{t=1,\dots,n_t}$, $h := (h_t)_{t=1,\dots,n_t}$, $\eta := (\eta_t)_{t=1,\dots,n_t}$, $\theta := (\theta_t)_{t=1,\dots,n_t}$, and $\alpha := (\alpha_t)_{t=1,\dots,n_t}$. The binary decision variables varying with time step $t$ are defined as: $v^+ := (v^+_t)_{t=1,\dots,n_t}$ and $v^- := (v^-_t)_{t=1,\dots,n_t}$. Finally, while the maximization of SCC can involve the installation of both PCVs and AFVs, flushing is not considered as an operational activity for pressure management. As a result, the formulation of valve placement and control problems for the minimization of AZP is generally restricted to the placement and control of PCVs (i.e. $y = \mathbb{0}^{n_n}$).

Following the notation in \citet{ECK2013}, we denote Problem \eqref{eq:minlp_problem} as VP-MINLP, which comprises both binary and continuous variables. Its solution represents the optimal placement and operational setting of PCV and AFV (if applicable) actuators. Moreover, we denote VC-NLP as the control-only subproblem of VP-MINLP. Here, decision variables include only optimal operational settings, while binary variables $v,\ y$ and $z$ are fixed to a given valve placement configuration, making Problem \eqref{eq:minlp_problem} a nonlinear programming (NLP) problem. VP-MINLP and VC-NLP are nonconvex optimization problems, where sources of nonconvexity include the frictional head loss model $\phi(\cdot)$ in \eqref{eq:hyd_energy_nonlin} and, for the self-cleaning problem, the sigmoidal SCC objective function in \eqref{eq:scc_objective_sigmoid}.

\subsection{Single-objective solution methods}
\label{sec:review_solution}
Previous literature has solved VP-MINLP using both metaheuristic and mathematical optimization methods. Metaheuristic approaches have included genetic algorithms \citep[e.g.][]{REIS1997,ARAUJO2006,CREACO2015,COVELLI2016}, harmony search \citep[e.g.][]{DEPAOLA2017}, and particle swarm optimization \citep[e.g.][]{DINI2020}. In order to comply with hydraulic constraints, these methods employ an external hydraulic solver (e.g. EPANET; \citeauthor{EPANET2.2}, \citeyear{EPANET2.2}) to evaluate the computed valve locations and settings. Mathematical optimization methods, on the other hand, embed the hydraulic equations as constraints within the optimization problem. These methods have included local optimization algorithms \citep{DAI2014,PECCI2021,PECCI2022a,PECCI2022b,JENKS2023} as well as global solvers \citep{ECK2013,DAMBROSIO2015,GHADDAR2017,PECCI2019}, where the latter provide certified optimality bounds for the computed feasible solutions. 

\begin{figure}[ht]
    \centering
    \includegraphics[width=0.6\textwidth]{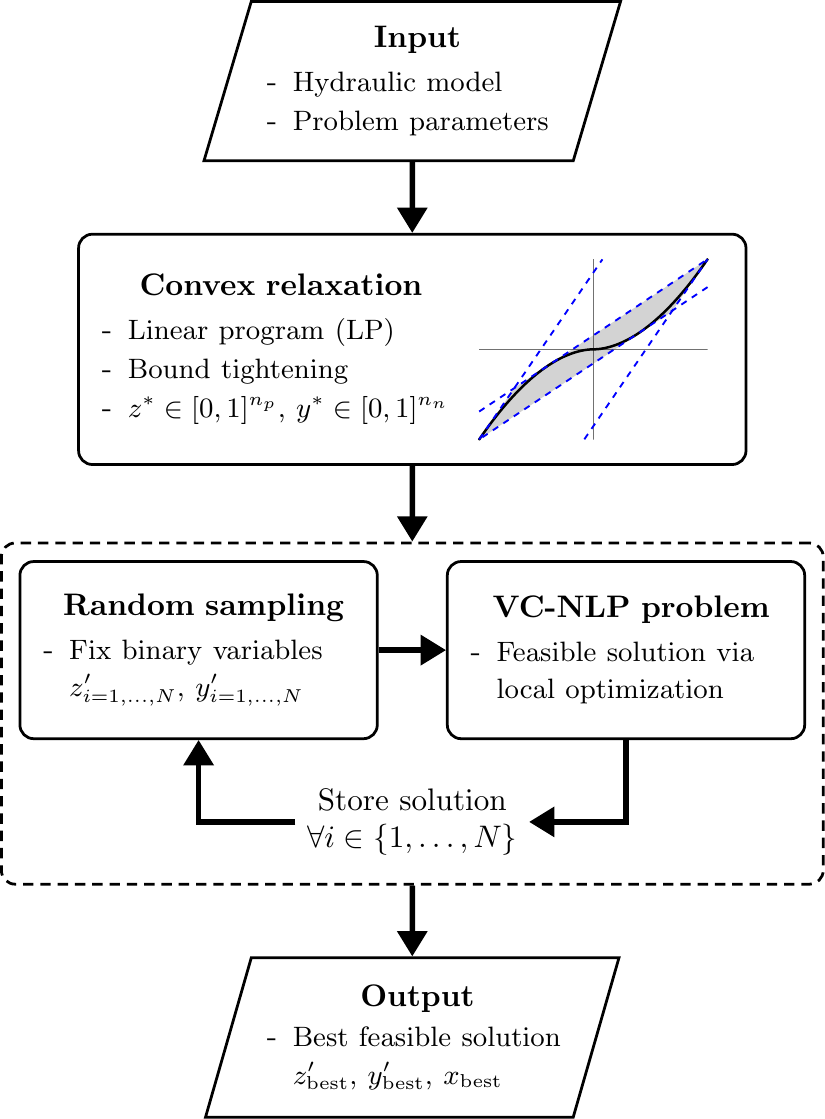}
    \caption{Convex heuristic algorithm.}
    \label{fig:convex_alg}
\end{figure}

This review focuses on heuristic methods developed to compute feasible solutions to large-scale MINLP problems. In particular, we consider the solution algorithm proposed in \citet{JENKS2023} for solving the SCC design-for-control problem. This algorithm comprised convex relaxations, a randomization heuristic, and a multi-start NLP solver to compute feasible solutions to VP-MINLP. Moreover, an optimization-based bound tightening (OBBT) scheme \citep{BELOTTI2009} was implemented to strengthen the convex relaxations. The algorithm extended the convex heuristic implemented in \citet{PECCI2022a} to include the nonconvex SCC objective function \eqref{eq:scc_objective} and AFVs as control actuators. \Cref{fig:convex_alg} provides an overview of the solution algorithm. The reader is referred to \citet{JENKS2023} for further details on the algorithm's implementation. 

Furthermore, we note that few studies have considered bi-directional (dynamic) boundary valves as control actuators in the optimization problem, which is essential for the design and operation of dynamically adaptive networks. Following \citet{JENKS2023}, Problem \eqref{eq:minlp_problem} permits bi-directional operation of control valve actuators through the set of binary variables $v_t \in \{0,1\}^{2n_p}$. As defined in \Cref{sec:review_hydraulics}, $v_t$ can take either positive or negative control direction for each time step $t \in \{1, \dots, n_t\}$, $v_t$. An enumeration approach was implemented to test all control direction combinations for a given valve configuration in the VC-NLP subproblem. That is, $v_t$ was fixed and VC-NLP was solved for $M = 2^{n_v}$ control direction combinations. Alternatively, \citet{ULUSOY2022a} formulated VP-MINLP with binary variables $v \in \{0,1\}^{2n_p}$ fixed in time $t$. A bi-linear constraint was then required to constrain $q_t$ and $\eta_t$ to the same direction at control valve links. Moreover, \citet{ULUSOY2022a} relaxed the resulting complementarity conditions on $q_t$ and $\eta_t$ by introducing an auxiliary variable. VP-MINLP and VC-NLP were then reformulated with a parametrized penalty term in the objective function.

For a feasible valve placement configuration $(v,y,z)$, the solution to VC-NLP provides a feasible (upper bound) solution of the aforementioned VP-MINLP. In addition, its solution is central to control scheduling problems where valve locations are known \textit{a priori} and only operational settings are optimized. In recent literature, two local solvers have been implemented to compute feasible solutions to VC-NLP: the state-of-the-art nonlinear solver IPOPT \citep[e.g.][]{ECK2013,DAI2014,GHADDAR2017,PECCI2019,NERANTZIS2020,ULUSOY2020}; and a sequential convex programming (SCP) algorithm \citep[e.g.][]{WRIGHT2014,WRIGHT2015,ABRAHAM2016,JENKS2023}. Application of these NLP solution algorithms have varied by problem formulation. For example, second-order derivatives, which are required for IPOPT, may not be well defined for some decision variables \citep[e.g. roughness coefficients for model calibration in][]{WALDRON2020}. In these cases, a tailored SCP algorithm has been used. The SCP algorithm comprises gradient-based optimization and a trust-region or line search strategy to evaluate solution progress. Furthermore, IPOPT can readily handle additional problem constraints, whereas SCP requires a tailored approach for each problem formulation to move from the current point $x_k$ to a new iterate $x_{k+1}$ in the solution process.

\subsection{Opportunities and challenges}
\label{sec:review_gaps}
The problem of integrating pressure management and self-cleaning controls has not been extensively investigated in the literature. This problem is particularly relevant as AZP and SCC have been identified as conflicting operational objectives \citep{ABRAHAM2016,ABRAHAM2019,JENKS2023}. Consequently, we have identified the following research opportunities and challenges, which we study in this paper.
\begin{itemize}
    \item \textbf{Bi-objective valve placement and control problem}. While AZP and SCC have been identified as conflicting objectives, an in-depth analysis of their trade-offs has not been performed for the valve placement and control problem. An investigation should assess the degree to which AZP and SCC are conflicting, especially when AFVs are considered as control actuators for the SCC problem. Moreover, the Pareto front should be approximated to evaluate the compromises between operational objectives.
    \item \textbf{Integrated control scheduling problem}. With valve actuators placed, the problem of integrating AZP and SCC controls within a control scheduling programme needs to be addressed. Since self-cleaning is a periodic mode of operation, optimal SCC controls may be implemented over a specific (and known) period. Otherwise, AZP can be considered the main (or continuous) operational objective.
\end{itemize}

The remainder of this paper is structured as follows. First, we formulate the bi-objective valve placement and control problem in \Cref{sec:bi_objective_formulation}. We apply the weighted sum scalarization method to approximate the Pareto optimal set (or Pareto front) between AZP and SCC objectives. In \Cref{sec:numerical_experiments}, we perform numerical experiments to demonstrate the AZP-SCC trade-offs for two case networks. Finally, we present a preliminary investigation of an adaptive control scheme for integrating AZP and SCC control modes.

\section{Bi-objective valve placement and control}
\label{sec:bi_objective_formulation}
This section formulates a bi-objective optimization problem to investigate the trade-offs between AZP and SCC. In doing so, we aim to search for optimal compromises between operational objectives, thereby informing decision making for the design and control of dynamically adaptive WDNs. This bi-objective problem extends the problem initially suggested in \citet{ABRAHAM2019} by introducing AFV actuator placement and control as decision variables. Note that AFVs are intended to replace (or complement) manual flushing practices. As a result, unlike pressure induced leakage, we do not consider the water discarded from the operation of AFVs as a form of undesirable water loss.

Using the same variable definitions in \eqref{eq:minlp_problem}, we formulate the bi-objective problem as 
\begin{linenomath}
\begin{subequations} \label{eq:minlp_bi_problem}
    \begin{alignat}{3}
    & \underset{\substack{x, \, v, \, z, \, y}}{\text{minimize}}
    & \quad & f(x) := (f_1(x),f_2(x)) \label{eq:minlp_bi_problem_a} \\
    & \text{subject to}
    & & Ax - \phi(x) = 0 \label{eq:minlp_bi_problem_b} \\
    & & & x \in X(v,y), \; v \in V(z)\label{eq:minlp_bi_problem_c} \\
    & & & v \in \{0,1\}^{2n_pn_t}, \; z \in \{0,1\}^{n_p}, \; y \in \{0,1\}^{n_n} 
\end{alignat}
\end{subequations}
\end{linenomath}
where $f_1$ and $f_2$ represent either AZP or $-$SCC. Following the notation in \Cref{sec:review_minlp}, we denote Problem \eqref{eq:minlp_bi_problem} as VP-BOMINLP and the bi-objective problem for a fixed valve configuration as VC-BONLP.

In this work, we implement the weighted sum (WS) scalarization method to approximate the Pareto front (or non-dominated set) of VP-BOMINLP. This is one of the most widely used scalarization method in the multi-objective optimization literature \citep{MIETTINEN1998,EHRGOTT2006}. However, a known limitation of the WS scalarization method is its inability to generate points in the nonconvex parts of the Pareto front \citep{DAS1997}. To this end, future work might investigate other scalarization methods to convert the bi-objective problem into a series of parametrized single-objective problems --- e.g. Chebyshev or $\epsilon$-constraint schemes \citep{EHRGOTT2006}. We reformulate VP-BOMINLP as a sequence of parametrized single-objective problems. The solutions of these single-objective problems define the set of Pareto optima of VP-BOMINLP, with its image in the objective space representing the approximated Pareto front. Since VP-BOMINLP is nonconvex, the Pareto front is defined only by local optima to the single-objective VP-MINLP problems. Here, such local optima to the single-objective problems are computed using the convex heuristic described in \Cref{sec:review_solution}. 

Let $x^*_i$ be the set of continuous variables that minimizes objective $f_i$. It follows that $f(x^*_i)$, $i=1,2$, denote the anchor points on the Pareto front. Moreover, $f_u := [f_1(x^*_1),f_2(x^*_2)]^T$ is the utopia point, which has values that simultaneously minimize both objective functions; by definition, $f_u$ is not a feasible solution to VP-BOMINLP. The objective functions are normalized using the utopia point and respective anchor points, as follows:
\begin{linenomath}
\begin{subequations}
    \label{eq:ws_normalization}
    \begin{align}
        &\bar{f}_1(x) := \frac{f_1(x) - f_1(x^*_1)}{f_1(x^*_2) - f_1(x^*_1)}  \label{eq:ws_normalization_a}\\
        &\bar{f}_2(x) := \frac{f_2(x) - f_2(x^*_2)}{f_2(x^*_1) - f_2(x^*_2)}  \label{eq:ws_normalization_b}
    \end{align}
\end{subequations}
\end{linenomath}
For $\omega \in [0,1]$, we then consider the WS problem formulation ${\text{VP-MINLP}_{\text{WS}}(\omega)}$:
\begin{linenomath}
\begin{subequations} \label{eq:minlp_ws_problem}
    \begin{alignat}{3}
    & \underset{\substack{x, \, v, \, z, \, y}}{\text{minimize}}
    & \quad & (1-\omega)\bar{f}_1(x) + \omega\bar{f}_2(x) \label{eq:minlp_ws_problem_a} \\
    & \text{subject to}
    & & Ax - \phi(x) = 0 \label{eq:minlp_ws_problem_b} \\
    & & & x \in X(v,y), \; v \in V(z)\label{eq:minlp_ws_problem_c} \\
    & & & v \in \{0,1\}^{2n_pn_t}, \; z \in \{0,1\}^{n_p}, \; y \in \{0,1\}^{n_n}.
\end{alignat}
\end{subequations}
\end{linenomath}
The Pareto front is generated by solving a sequence of parametrized problems ${\text{VP-MINLP}_{\text{WS}}(\omega)}$ for a given set of weights $\omega$. In \Cref{sec:num_experiment_bi_dfc}, we compute points on the Pareto front for $n=10$ evenly distributed weights in the domain $[0,1]$.

\section{Numerical experiments}
\label{sec:numerical_experiments}
Numerical experiments reported herein were performed in Julia 1.8.2 \citep{BEZANSON2017} with optimization solvers accessed via JuMP 1.4.0 \citep{DUNNING2017}. All computations were performed within the Ubuntu 20.04.5 LTS Linux distribution on a 2.50-GHz Intel(R) Core(TM) i9-11900H CPU with 8 cores and 32.0 GB of RAM. Linear programs were solved using Gurobi 9.5.2 \citep{GUROBI2022}. Nonlinear programs were solved using IPOPT 3.14.4 \citep{WACHTER2006} installed with the HSL MA57 linear solver \citep{HSL2021}.

\begin{figure}[t]
    \centering
    \subfloat[\centering \label{fig:case_networks_modena}\texttt{Modena}]{{\includegraphics[width=0.43\textwidth]{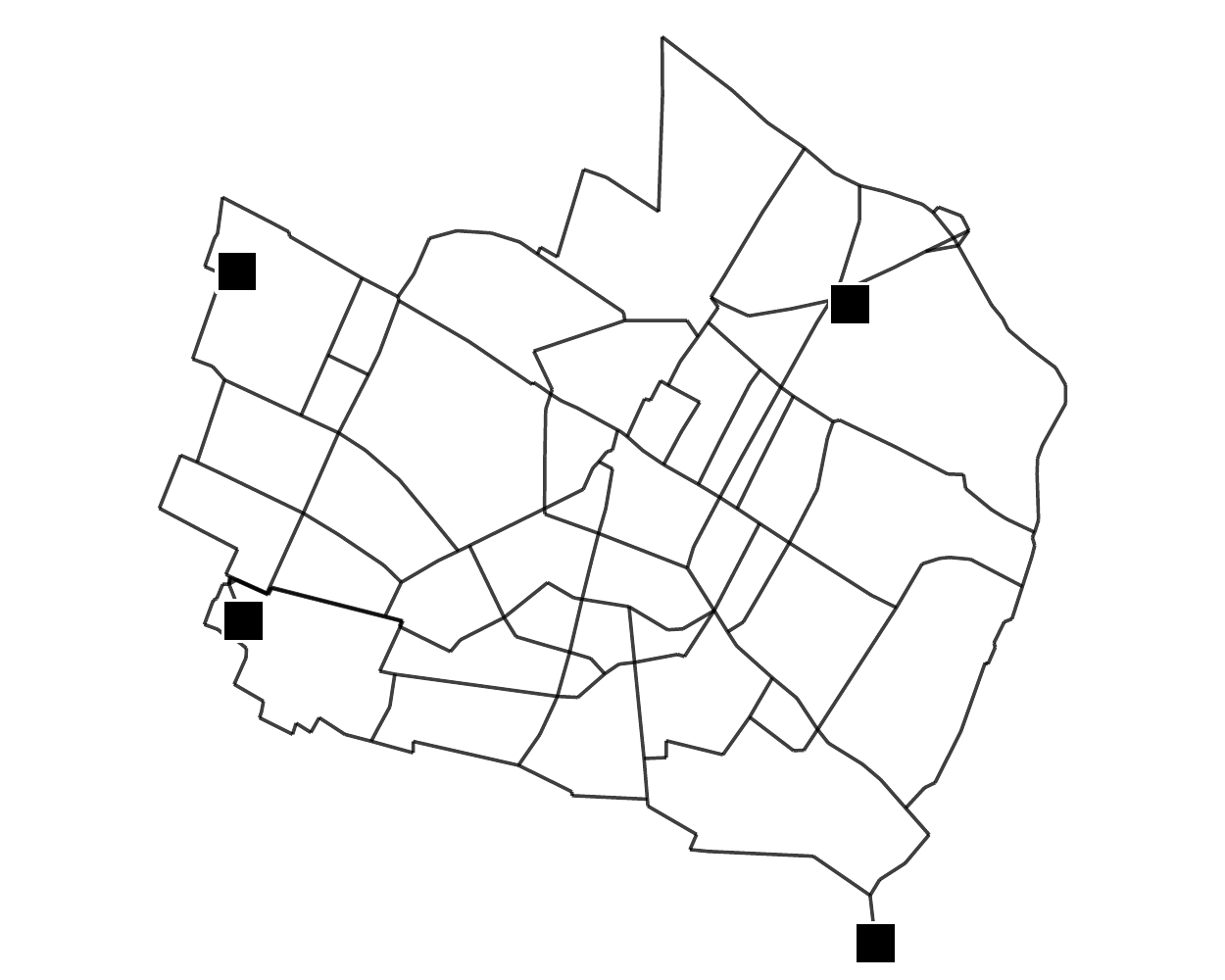} }}
    \;
    \subfloat[\centering \label{fig:case_networks_bwfl}\texttt{BWFLnet}]{{\includegraphics[width=0.53\textwidth]{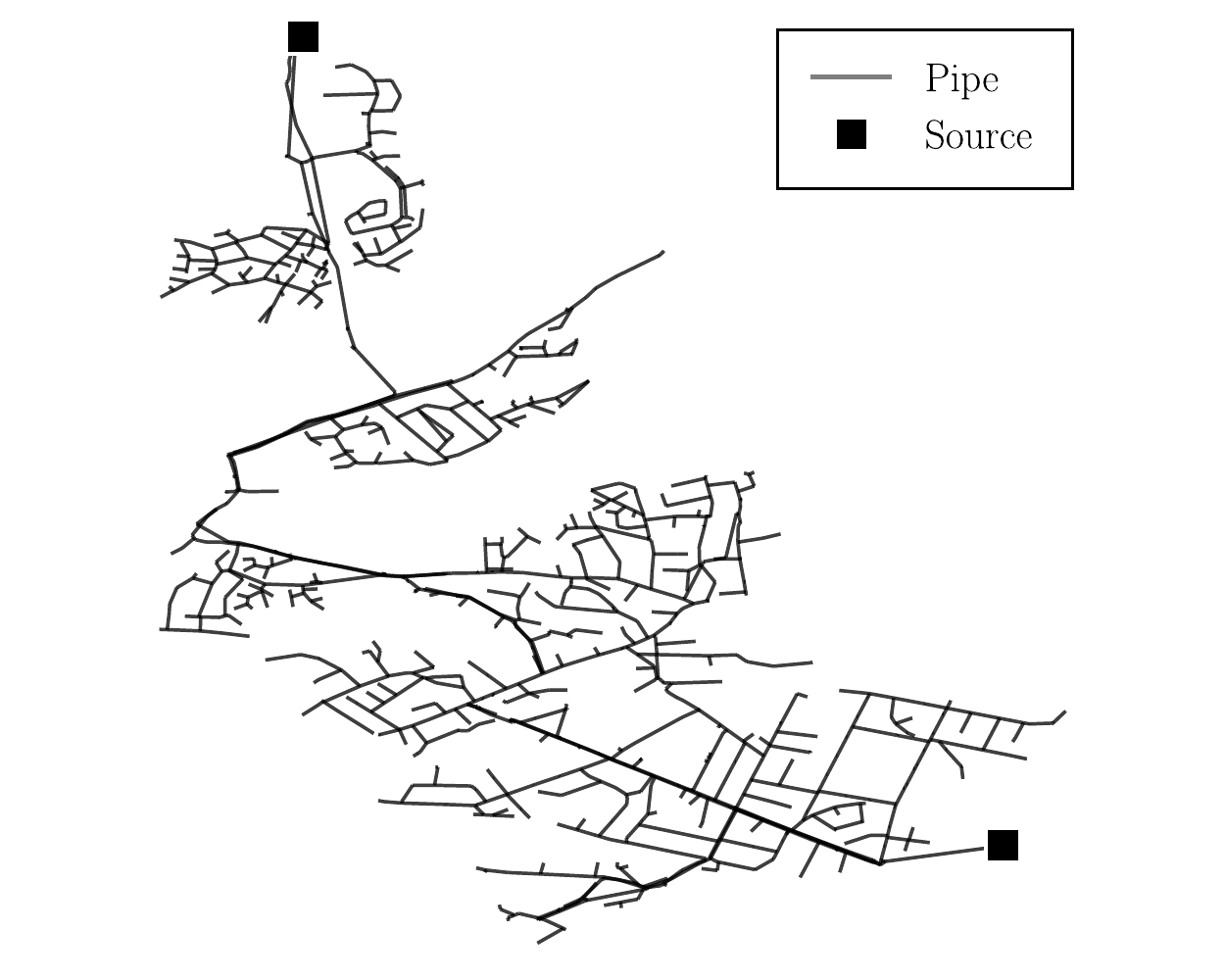} }}
    \caption{Case network layouts.}
    \label{fig:case_networks}
\end{figure}

We consider two case networks in this study: \texttt{Modena} and \texttt{BWFLnet}. \texttt{Modena} is a reduced version of the hydraulic model for a medium-sized Italian city and is a well-established benchmarking network in the literature \citep{BRAGALLI2012}. \texttt{Modena} consists of $n_n=268$ demand or junction nodes, $n_p=317$ links, and $n_0=4$ source nodes. Moreover, it is modelled over $n_t=24$ hourly discrete time steps. It has a highly interconnected (looped) network topology and large demands, representative of flows typically conveyed at the transmission level. Its layout is shown in \Cref{fig:case_networks_modena}. In comparison, \texttt{BWFLnet} is a large-scale operational network in the UK, which is jointly operated by Imperial College London, Bristol Water Plc, and Cla-Val Ltd~\citep{WRIGHT2014}. \texttt{BWFLnet} is a unique smart water network demonstrator, enabling the dynamically adaptive control of network topology and hydraulic conditions. It has been used as a case network in numerous optimization and control studies \citep[e.g.][]{WRIGHT2014,WRIGHT2015,PECCI2019,NERANTZIS2020,PECCI2022a,JENKS2023}. The \texttt{BWFLnet} network used in this study consists of $n_n=2,745$ demand or junction nodes, $n_p=2,816$ links, and $n_0=2$ source nodes. It is modelled over $n_t=96$ 15-minute discrete time steps. In contrast to \texttt{Modena}, \texttt{BWFLnet} represents an operational network with branched topology (i.e. little redundancy in connectivity). Moreover, its demands represent water consumption in a typical UK urban area as opposed to the aggregated demands consistent with the magnitude of flows at the transmission level. \textcolor{red}{}The \texttt{BWFLnet} layout is shown in \Cref{fig:case_networks_bwfl}. Together, these case networks represent varied network conditions to demonstrate the investigated AZP-SCC control challenges.

The following assumptions are adopted for both single- and bi-objective optimization problems. Frictional head losses are modelled using the HW equation and the known roughness coefficients from each network model. The minimum pressure head is set to 15 m at demand nodes and 0 m otherwise. Moreover, the maximum flushing rate is set to 25 L/s at AFV nodes, which corresponds to the recommended fire flow capacity at hydrants in the UK \citep{FIRE2007}. For the approximated SCC objective function in \eqref{eq:scc_objective_sigmoid}, we set the sigmoid curvature parameter $\rho = 50$ and the self-cleaning velocity threshold $u_j^{\min} = 0.2$ m/s for all $j \in \{1, \dots, n_p\}$. The selected velocity threshold $u_j^{\min}$ references experimental findings reported in the literature \citep{RYAN2008,BLOKKER2010}. No existing PCV or AFV actuators are considered for the single- and bi-objective design-for-control problems in \Cref{sec:num_experiment_single_dfc} and \Cref{sec:num_experiment_bi_dfc}, respectively. However, the adaptive control scheme investigated in \Cref{sec:num_experiment_adaptive_control} assumes a known valve configuration. Furthermore, we apply a modified version of the convex heuristic from \citet{JENKS2023} to compute feasible solutions to the single-objective VP-MINLP problem. Here, the state-of-the-art NLP solver IPOPT is used in place of the sequential convex programming algorithm. We note that, in contrast to that reported in \citet{JENKS2023}, no infeasibility errors were observed when using IPOPT 3.14.4 (accessed via JuMP) for the current study. Moreover, only one starting point is tested when solving VC-NLP, corresponding to the valve settings returned by the solution of the convex relaxation in VP-MINLP (see \Cref{fig:convex_alg}). Finally, hydraulic equations \eqref{eq:hyd_energy_nonlin}-\eqref{eq:hyd_mass} are solved using the null space solver developed by \citet{ABRAHAM2015}.

\subsection{Single-objective optimization}
\label{sec:num_experiment_single_dfc}
\begin{figure}[t]
    \centering
    \subfloat[\centering \label{fig:modena_azp_dfc}Minimization of AZP in \texttt{Modena}.]{
        \includegraphics[width=0.44\textwidth]{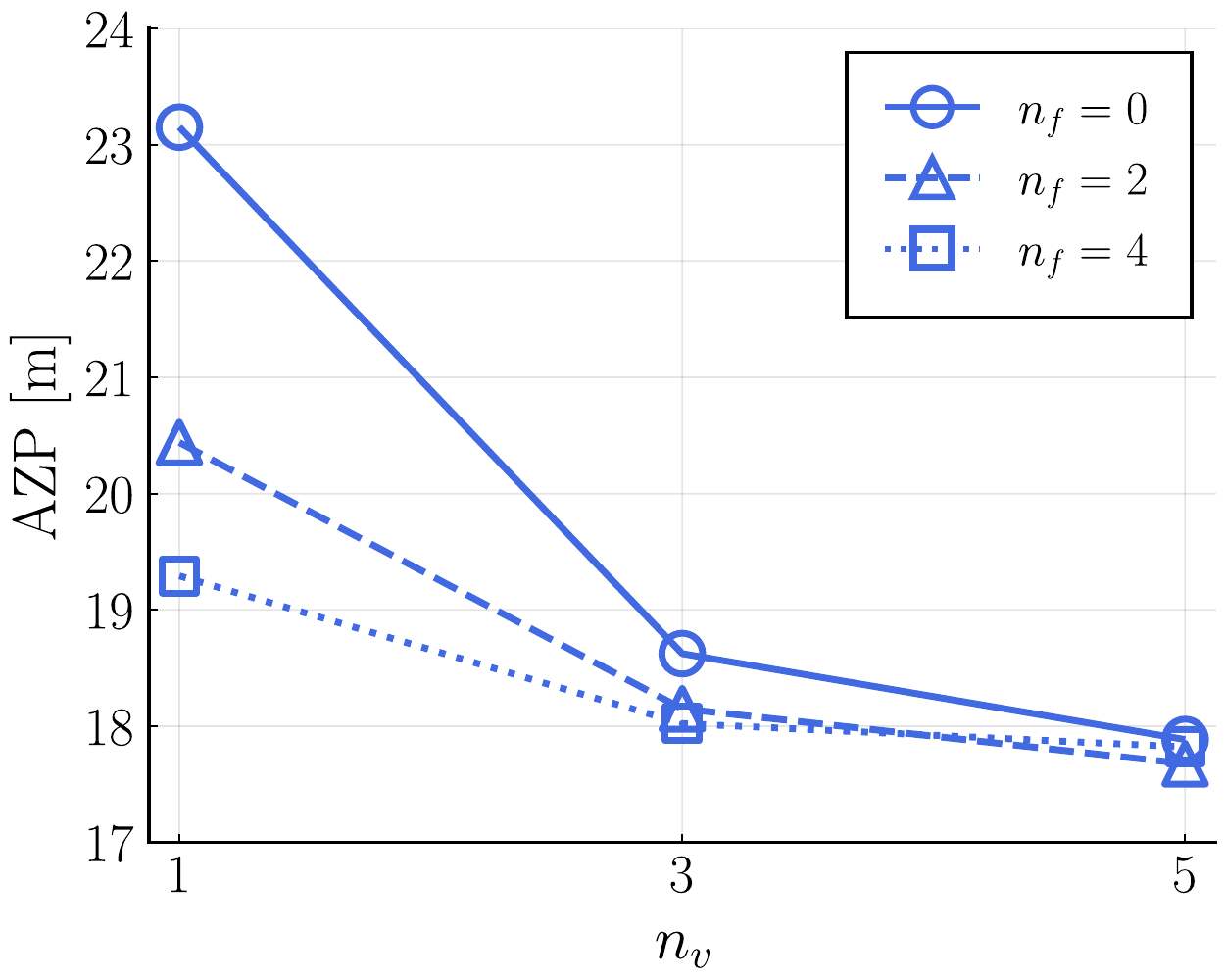}}
    \qquad
    \subfloat[\centering \label{fig:modena_scc_dfc}Maximization of SCC in \texttt{Modena}.]{
        \includegraphics[width=0.44\textwidth]{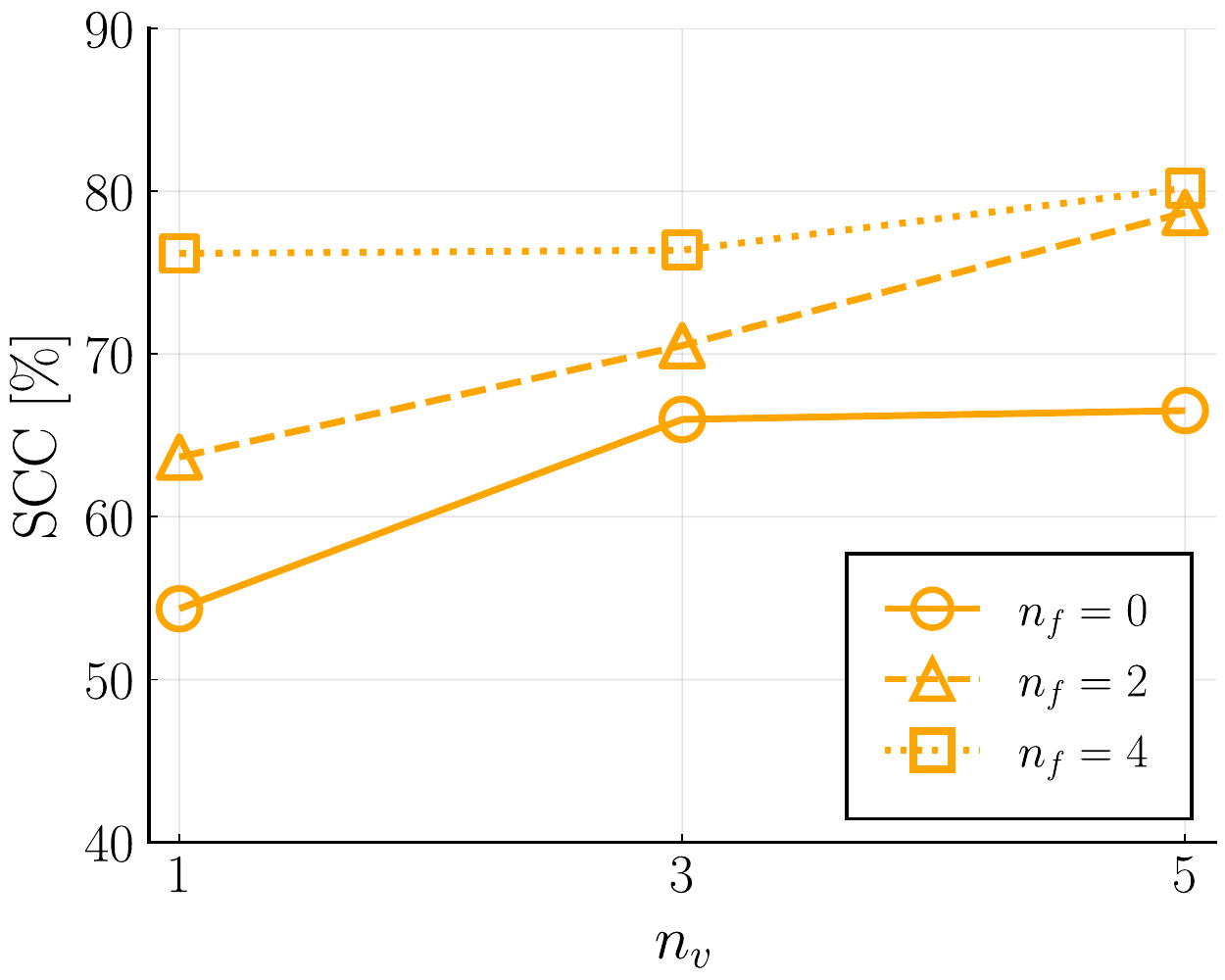}}
    \\ \vspace{0.2cm}
    \subfloat[\centering \label{fig:bwfl_azp_dfc}Minimization of AZP in \texttt{BWFLnet}.]{
        \includegraphics[width=0.44\textwidth]{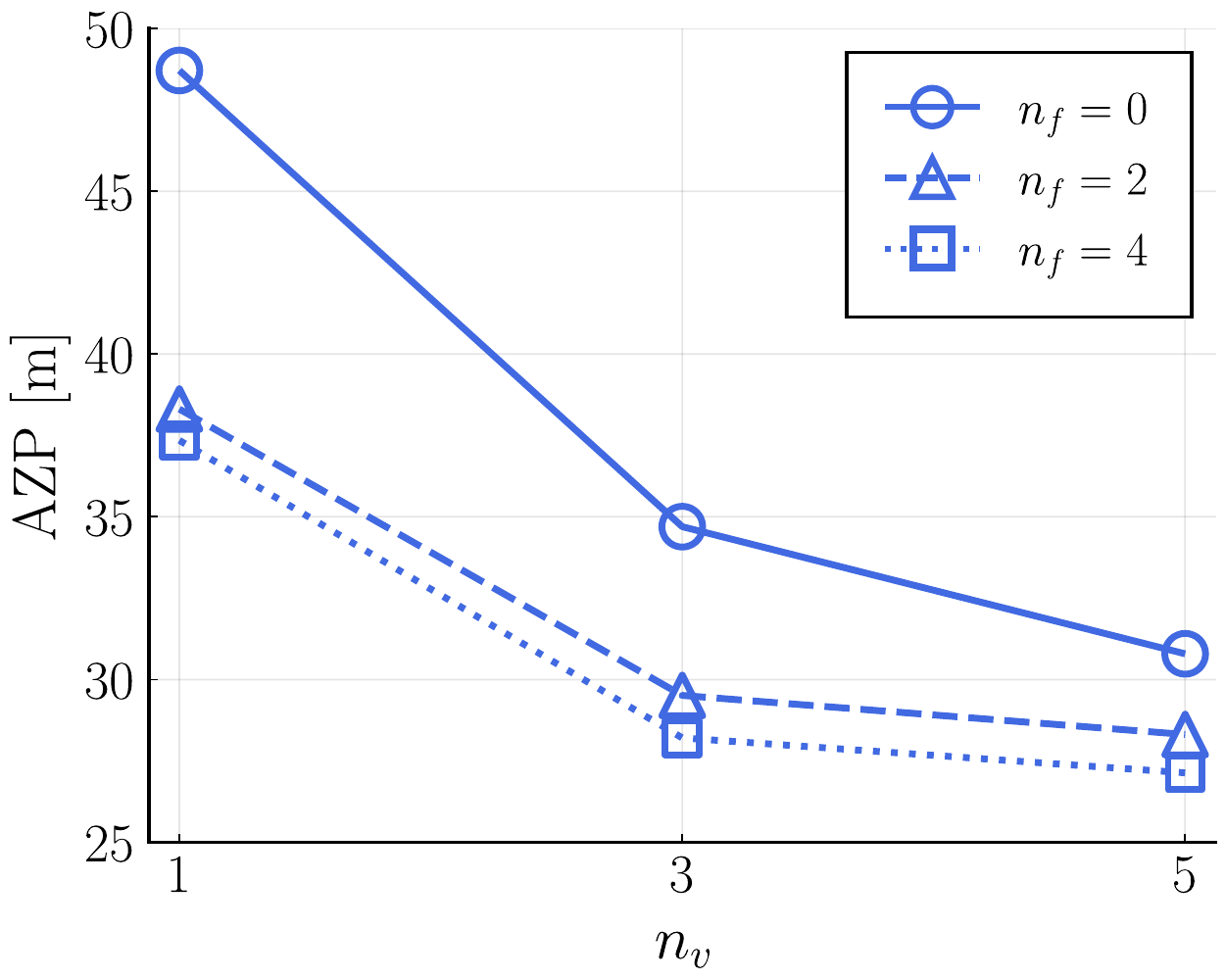}}
    \qquad
    \subfloat[\centering \label{fig:bwfl_scc_dfc}Maximization of SCC in \texttt{BWFLnet}.]{
        \includegraphics[width=0.44\textwidth]{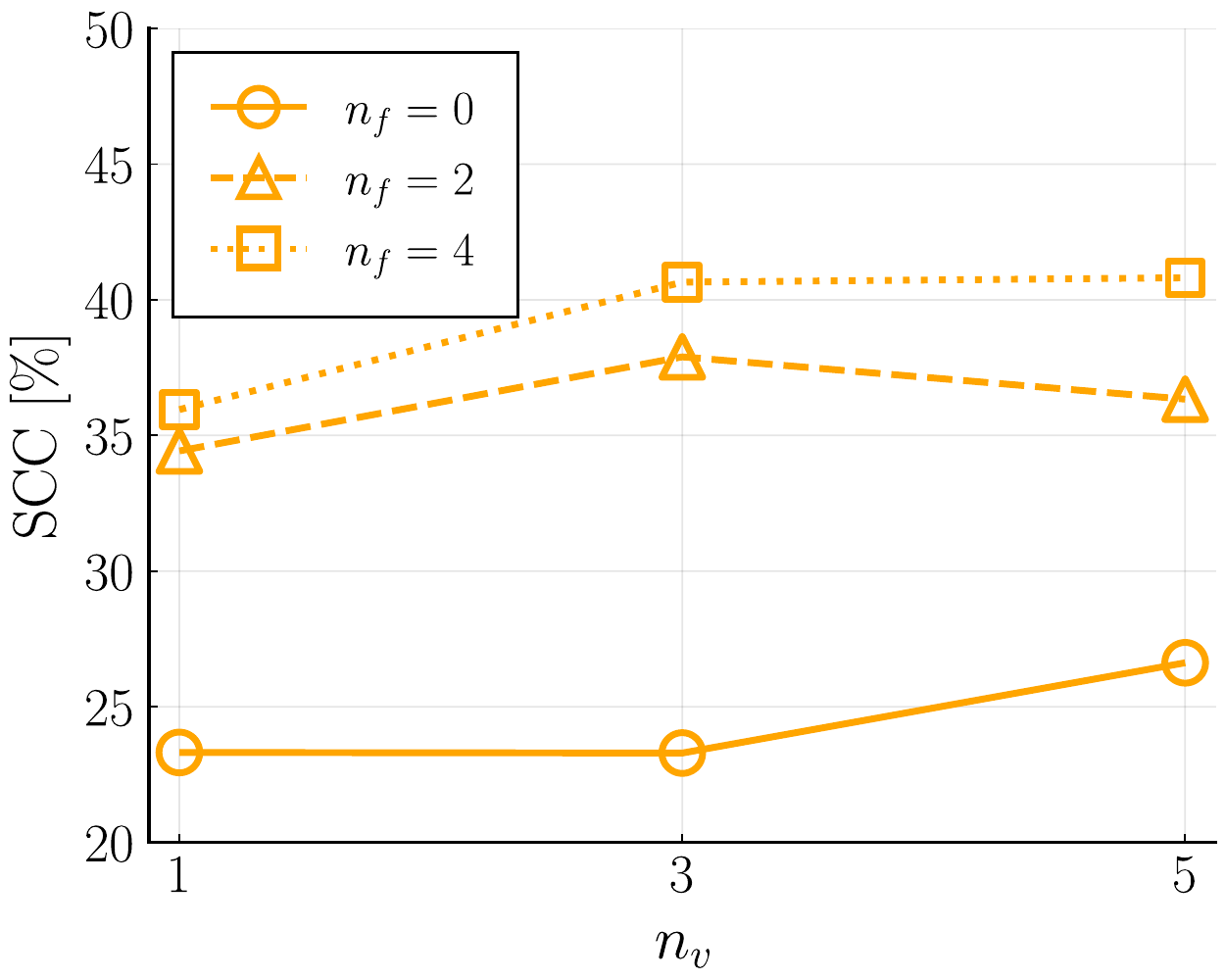}}
    \caption{Individual AZP and SCC single-objective VP-MINLP results.}
    \label{fig:single_obj_dfc}
\end{figure}

We first compute feasible solutions to VP-MINLP for individually optimizing AZP and SCC in \texttt{Modena} and \texttt{BWFLnet}. We consider the formulation of \eqref{eq:minlp_problem} for valve configurations comprising $n_v \in \{1, 3, 5\}$ PCV and $n_f \in \{0, 2, 4\}$ AFV actuators. These results seek to find the best number of $n_v$ and $n_f$ control actuators for consideration in the bi-objective problem formulation. Furthermore, the objective values of the solutions to the single-objective problems are used to normalize the objective functions in the weighted sum scalarization method. We note that, in practice, AFVs are not considered as part of a pressure management programme. However, they are included as valve actuators in the single-objective AZP problems to support the bi-objective problem formulation. This is primarily because the AZP-SCC bi-objective problem formulation in \Cref{sec:bi_objective_formulation} requires a pre-defined valve configuration, of which considers both PCV and AFV actuators. \Cref{fig:single_obj_dfc} presents the optimal valve placement and control results individually obtained for AZP and SCC in case networks \texttt{Modena} and \texttt{BWFLnet}.

We observe a reduction in AZP of approximately $23\%$ for \texttt{Modena} and $38\%$ for \texttt{BWFLnet} between problem instances with $n_v=1$ and $n_v=3$ PCVs. However, in both networks, additional benefits appear to diminish for the case with $n_v=5$. These results are generally in accordance with previous studies concerning optimal valve placement and control for minimizing AZP ~\citep{PECCI2019,PECCI2022a}. Accordingly, we assume the valve configuration corresponding to $n_v=3$ to represent the most economical design decision of the evaluated problem instances. The addition of AFV actuators had varying levels of influence on AZP performance. For cases with $n_v > 1$, flushing demands showed little-to-no improvements in \texttt{Modena}. On the other hand, a noticeable reduction in AZP is discerned for problem instances with AFVs in \texttt{BWFLnet}. This observation can be attributed to the additional head losses materializing from flushing demands, particularly in areas of the network where elevation or hydraulic distance inhibits pressure modulation from PCVs. As noted previously, we include the AZP problem results with AFV actuators to support the AZP-SCC bi-objective problem formulation. 

Results for the single-objective SCC problem indicate that the optimal placement and control of AFV actuators deliver the greatest benefits for improving self-cleaning performance. In the absence of flushing, SCC improvements are achieved via changes in network connectivity and pressure modulation at PCVs. Such dynamic changes in connectivity enable the redistribution of flow to maximize self-cleaning velocities. This is most evident in the highly interconnected (looped) topology of \texttt{Modena}, where problem instances of $n_v > 1$ show a noticeable improvement (roughly $20\%$) in SCC performance. On the other hand, the branched topology of \texttt{BWFLnet} limits control on the distribution of flow, as demonstrated by the relatively small increase in SCC performance for the problem instance of $n_v=5$. Moreover, we generally observe an increasing relationship between SCC performance and the number of AFVs considered in the problem formulation ($n_f$). This is to be expected for highly branched networks as the influence of flushing demands is limited to the upstream flow path.

\begin{figure}[p!]
    \centering
    \subfloat[\centering \label{fig:modena_net_dfc}\texttt{Modena} ]{\includegraphics[width=0.55\textwidth]{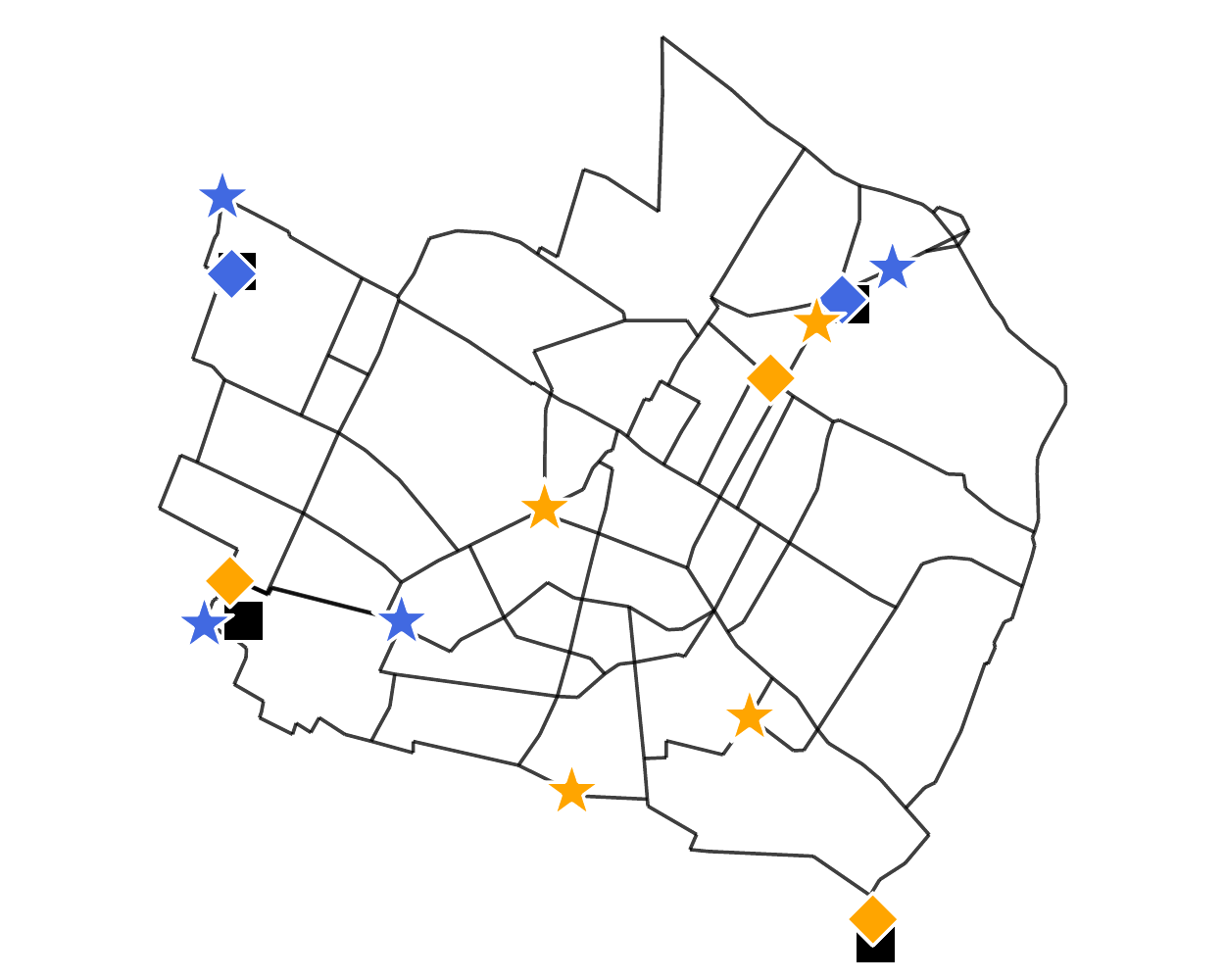} }
    \\ \vspace{1.0cm}
    \subfloat[\centering \label{fig:bwfl_net_dfc}\texttt{BWFLnet}]{{\includegraphics[width=0.62\textwidth]{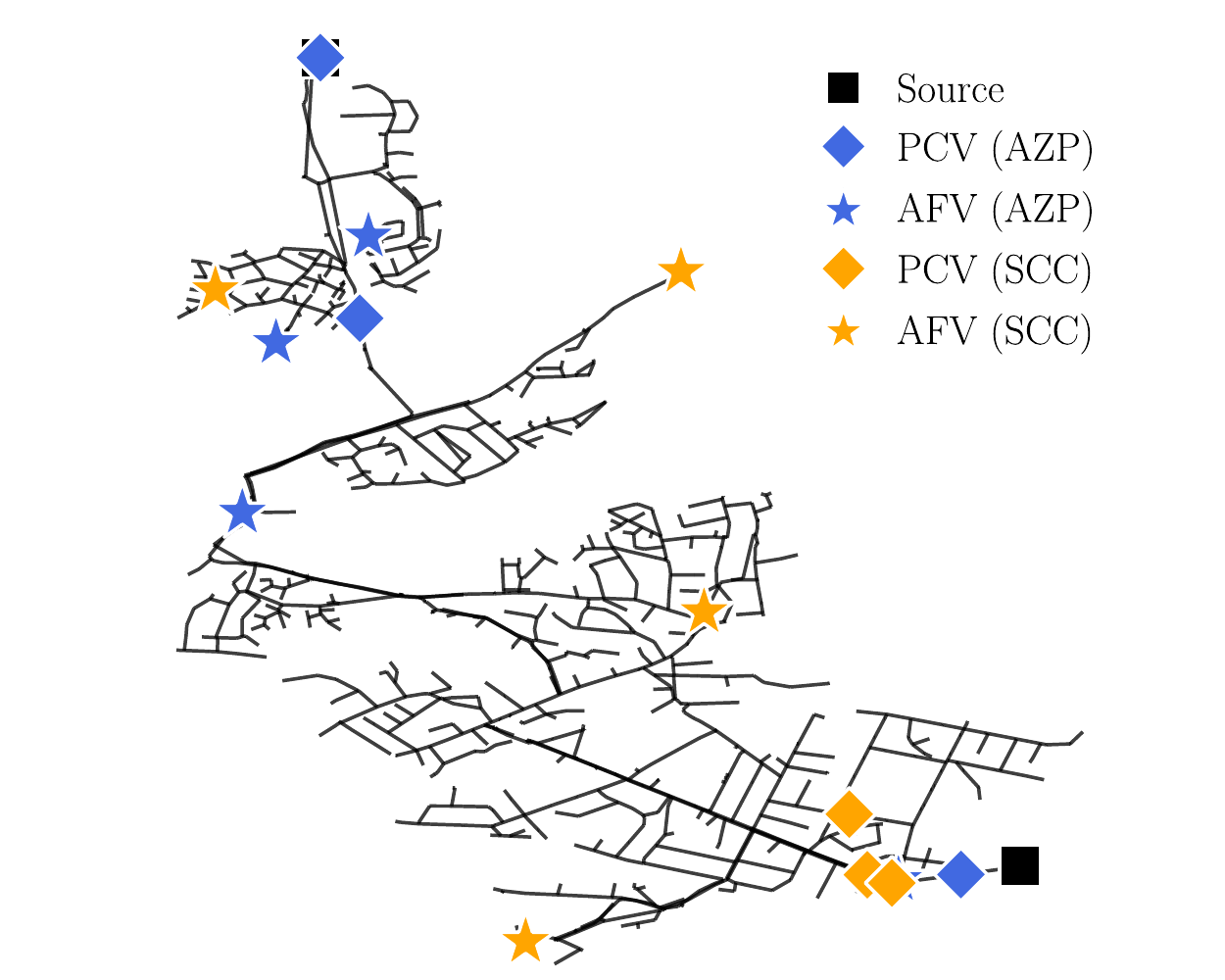} }}
    \vspace{0.5cm}
    \caption{Optimal pressure control (PCV) and automatic flushing (AFV) valve actuator locations for single-objective \textcolor{royal_blue}{\textbf{AZP}} and \textcolor{orange_yellow}{\textbf{SCC}} VP-MINLP problems.}
    \label{fig:single_obj_design}
\end{figure}

On the basis of these results, we consider the valve configuration of $n_v=3$ and $n_f=4$ for the bi-objective problems investigated in \Cref{sec:num_experiment_bi_dfc}. \Cref{table:azp_scc_anchor} compares the AZP and SCC objective values resulting from their respective optimal valve configurations (i.e. anchor points for the bi-objective problem formulation).
\begin{table}[t]
\centering
\setlength{\tabcolsep}{6pt}
\caption{Comparison of AZP and SCC anchor point objective values.}
\label{table:azp_scc_anchor}
    \begin{tabular}{lrr}
        \toprule
        Bi-objective anchor point & AZP [m] & SCC [\%] \\
        \midrule
        \multicolumn{3}{@{}l}{\texttt{Modena}} \\
        AZP problem & 18.0 & 59.2 \\
        SCC problem & 29.3 & 78.4 \\
        \midrule
        \multicolumn{3}{@{}l}{\texttt{BWFLnet}} \\
        AZP problem & 28.2 & 28.5 \\
        SCC problem & 48.9 & 42.0 \\
        \bottomrule
    \end{tabular}
\end{table} 
As expected, these objective values show significant differences between anchor points for both \texttt{Modena} and \texttt{BWFLnet}. This confirms the conclusions reported in \citet{ABRAHAM2016,ABRAHAM2019,JENKS2023} that AZP and SCC are indeed conflicting objectives. However, an analysis of their trade-offs is needed to determine if an acceptable compromise between objectives can be achieved. To gain a better understanding of these conflicts, we first review the corresponding optimal valve locations for the individual AZP and SCC problems, which are shown in \Cref{fig:single_obj_design}. We observe optimal PCV locations for both AZP and SCC problems to be near the source nodes (e.g. reservoirs) in the looped \texttt{Modena} network. In contrast, PCV locations are quite different in the more branched \texttt{BWFLnet} network, which show a cluster of PCVs for the SCC problem. This might suggest less marked trade-offs between AZP and SCC in \texttt{Modena} compared to \texttt{BWFLnet}. On the other hand, no clear pattern is discerned for optimal AFV locations. Future work is needed to further investigate the influence of network topology and valve placement.

\subsection{Bi-objective optimization}
\label{sec:num_experiment_bi_dfc}
In this section, we investigate the trade-offs between AZP and SCC objectives through the notion of Pareto optimality. This is intended to inform design decisions for the placement and operation of PCV and AFV actuators. In particular, we propose a hierarchical design strategy that leverages the locations of optimal valve actuators from the single-objective AZP and SCC problems presented in \Cref{sec:num_experiment_single_dfc}. For the AZP-SCC bi-objective problem formulated in \Cref{sec:bi_objective_formulation}, we implement the weighted sum (WS) scalarization method to compute the set of Pareto optima for $n=10$ evenly distributed weights in the domain $[0,1]$. Feasible solutions to the resulting sequence of single-objective optimization problems are computed using the convex heuristic described in \Cref{sec:review_solution}. Since~\eqref{eq:minlp_bi_problem} is nonconvex, the WS scheme generates an approximation of the Pareto front to the bi-objective problem. We discuss the benefits and limitations of the solution method in the subsequent results, particularly in relation to the nonconvexity of the considered problems. 

We compute approximated Pareto fronts for two bi-objective problem cases. The main case considers a valve control bi-objective NLP (VC-BONLP) problem, which has a known valve configuration corresponding to a hierarchical design strategy. We formulate VC-BONLP by fixing valve actuator locations using the single-objective results presented in \Cref{sec:num_experiment_single_dfc}. More specifically, we set the locations of PCVs based on the solution of the single-objective problem minimizing AZP. Then, with PCVs already configured, we set the locations of AFVs based on the solution of the single-objective problem maximizing SCC. Since the single-objective results demonstrate that AZP and SCC depend strongly on PCV and AFV controls, respectively, we test the hypothesis that this hierarchical design strategy can yield sufficient performance for both objectives. Additionally, we consider the design-for-control bi-objective problem (VP-BOMINLP), where valve locations and operational settings are jointly optimized. The solution to VP-BOMINLP aims to qualify the results obtained for VC-BONLP and demonstrate the trade-offs between AZP and SCC for different valve configurations.  

\begin{figure}[h!]
    \centering
    \subfloat[\centering \label{fig:modena_pareto}\texttt{Modena} ]{
    \includegraphics[width=0.44\textwidth]{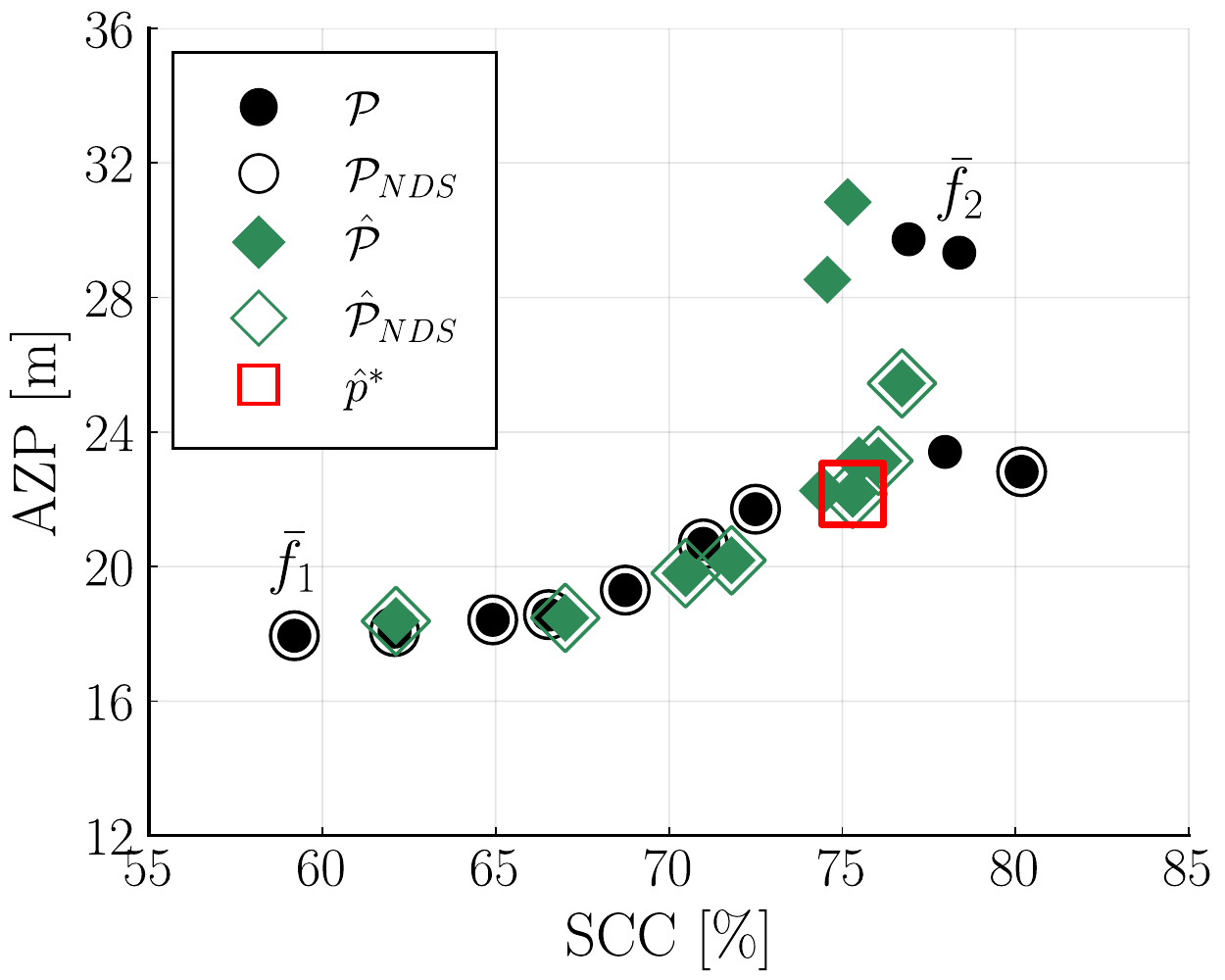}}
    \qquad
    \subfloat[\centering \label{fig:bwfl_pareto}\texttt{BWFLnet} ]{
    \includegraphics[width=0.44\textwidth]{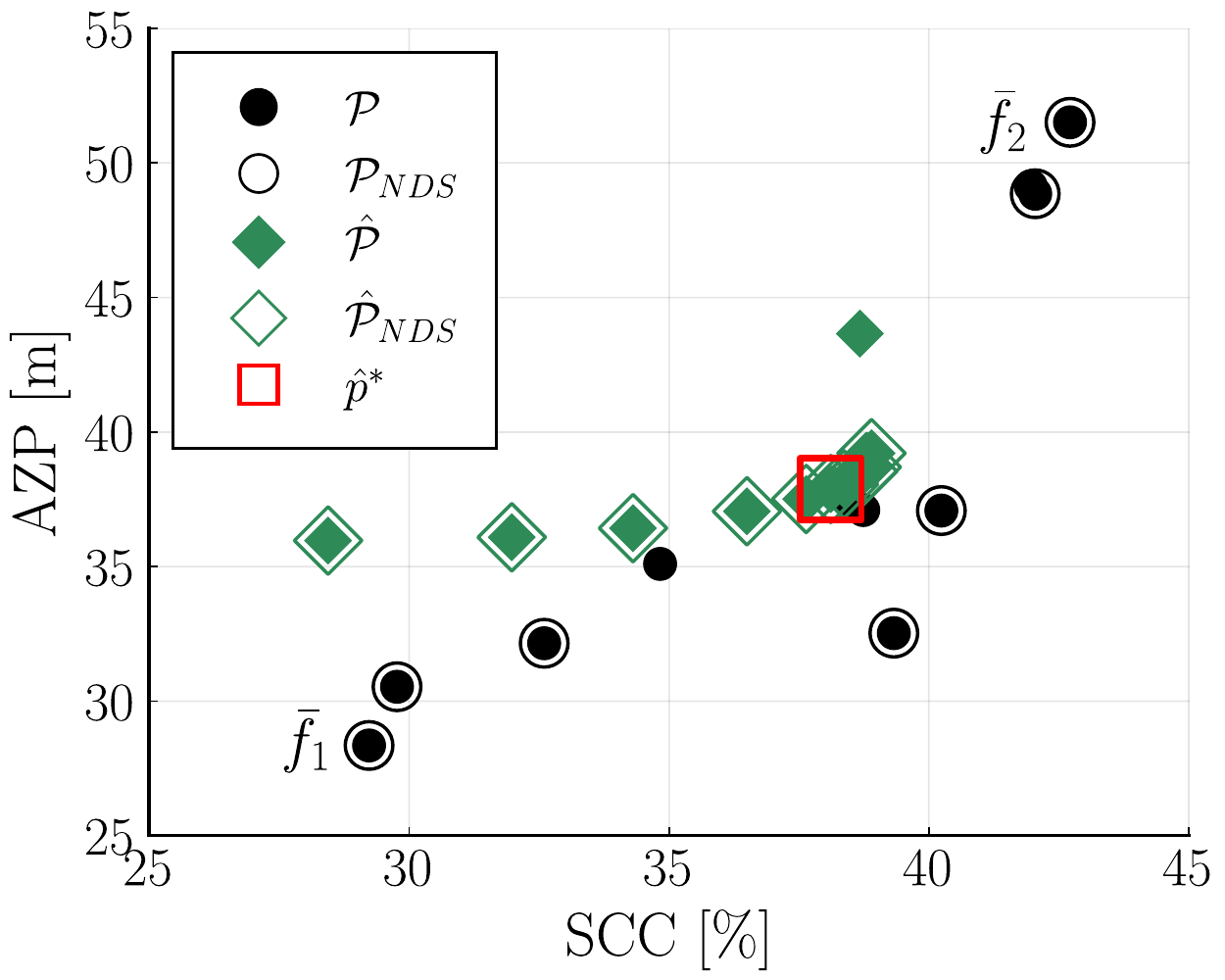}}
    \caption{Solutions $\mathcal{P}$ and $\hat{\mathcal{P}} \subseteq \mathcal{P}$ computed, respectively, for VP-BOMINLP and VC-BONLP with $n_v=3$ and $n_f=4$ using WS. Non-dominated solution sets $\mathcal{P}_{NDS}$ and $\hat{\mathcal{P}}_{NDS}$ are highlighted with offset borders and anchor points $\bar{f}_1$ and $\bar{f}_2$ are labelled for $\mathcal{P}_{NDS}$.}
    \label{fig:pareto}
\end{figure}

The following notation is used for presenting the bi-objective results. We denote the set of solutions computed by the WS method for VP-BOMINLP and VC-BONLP as $\mathcal{P}$ and $\hat{\mathcal{P}}$, respectively. Using the Pareto filter proposed in \citet{MESSAC2003}, we then identify the non-dominated subset of these solutions, which we denote as $\mathcal{P}_{NDS} \subseteq \mathcal{P}$ and $\hat{\mathcal{P}}_{NDS} \subseteq \hat{\mathcal{P}}$. These non-dominated sets approximate the Pareto front of the respective AZP-SCC bi-objective problems. The resulting Pareto fronts approximated for VP-BOMINLP and VC-BONLP problems with $n_v=3$ and $n_f=4$ are shown in \Cref{fig:pareto}. Observe that the Pareto fronts are fairly uniform and well-spread, suggesting that the WS scalarization method and selected weights are sufficient for evaluating trade-offs between objectives.

\Cref{fig:pareto} illustrates the trade-offs between AZP and SCC for the evaluated networks. For both \texttt{Modena} and \texttt{BWFLnet}, we observe that the sets $\mathcal{P}_{NDS}$ and $\hat{\mathcal{P}}_{NDS}$ describe relatively flat Pareto fronts progressing from the AZP anchor point in the bottom left corner, before incurring a noticeable increase in AZP in order to achieve conditions near the SCC anchor point. The selection of solution $\hat{p}^* \in \hat{\mathcal{P}}_{NDS}$, denoted by the red square in \Cref{fig:bwfl_pareto} and \Cref{fig:modena_pareto}, highlights that significant improvements in SCC can mostly be achieved for only moderate increases in AZP. Moreover, this behaviour is more more pronounced in \texttt{BWFLnet}, where the zone of influence of PCV actuators can be limited by the branched topology of the network (as shown in \Cref{fig:bwfl_net_dfc}). The looped topology of \texttt{Modena}, on the other hand, results in a wider area of the network affected by pressure modulation. Therefore, the trade-off in AZP is less marked when achieving hydraulic conditions corresponding to the SCC anchor point. Furthermore, we observe that some solutions returned by WS for VP-BOMINLP and, to a lesser extent, VC-BONLP are largely dominated by others. This is especially prevalent for the solutions of scalarized problems corresponding to larger SCC weights, for which the likelihood of computing suboptimal local solutions is increased by the greater degree of nonconvexity of the SCC objective compared to AZP. This phenomenon is applicable to both VP-BOMINLP and VC-BONLP problems and could potentially be mitigated through a more robust solution methodology --- e.g. multi-start approach proposed in \citet{JENKS2023}. 

Finally, we compare Pareto fronts of VC-BONLP and VP-BOMINLP. As expected, we observe an acceptable range in objective function values when PCV and AFV actuator locations are fixed according to the aforementioned hierarchical design. This highlights the strong dependence of each valve actuator on their respective operational objectives. \Cref{fig:modena_pareto} shows a small optimality gap between VC-BONLP and VP-BOMINLP anchor points for \texttt{Modena}. In comparison, \Cref{fig:bwfl_pareto} reveals slightly larger deviations between approximated Pareto sets $\mathcal{P}_{NDS}$ and $\hat{\mathcal{P}}_{NDS}$ for \texttt{BWFLnet}. In particular, we observe an 8~m difference in AZP values between the anchor points in the bottom left corner. Although this signifies a noticeable optimality gap between VC-BONLP and VP-BOMINLP, we note that flushing is not considered as an operational activity for pressure management. Therefore, a more reasonable comparison of AZP performance is with the single-objective problem instance of $n_v=3$ and $n_f=0$ in \Cref{fig:bwfl_azp_dfc}, which yields similar performance to the set $\hat{\mathcal{P}}_{NDS}$. 

The generated Pareto fronts allow an operator to (i) visualize the conflicting nature of AZP and SCC objectives and (ii) quantify such trade-offs to inform decision making for designing WDNs for the simultaneous control of pressure and water quality. In particular, we observe that a hierarchical design strategy (Problem VC-BONLP) yields good trade-offs between AZP and SCC objectives by varying actuator settings for pre-defined valve locations. In this approach, PCVs are first optimally placed to minimize AZP. Then, using these PCV locations, AFVs are placed to augment SCC conditions in the network. Informed by the Pareto front of VC-BONLP, valve actuator settings can be selected to balance the trade-offs between objectives. This, however, considers a static operational point for integrating AZP and SCC controls. In contrast to this approach, previous literature has suggested the occurrence of a daily (or every other day) self-cleaning period to be sufficient for mitigating discolouration risk in WDNs \citep{VREEBURG2009,BLOKKER2012}. On this basis, we describe a control scheme that facilitates the dynamic transition between AZP and SCC controls in \Cref{sec:num_experiment_adaptive_control}.

\subsection{Adaptive control scheme}
\label{sec:num_experiment_adaptive_control}
This section investigates the second research problem identified in \Cref{sec:review_gaps}: an integrated control problem for reducing AZP and increasing SCC in operational networks. We formulate an adaptive control scheme as a sequence of single-objective problems, whereby control modes are switched at pre-defined periods. This follows a similar approach to the model predictive control (MPC) scheme proposed in \citet{NERANTZIS2022}, which sought to adapt network hydraulics for fire fighting activity. For integrating AZP-SCC controls, however, we assume known periods for switching between operational modes, thereby facilitating proactive network self-cleaning maintenance. Furthermore, since background leakage is a continuous phenomenon occurring throughout the network (see review in \Cref{sec:review_obj_pressure}), we set AZP as the primary (or master) control mode and SCC as the secondary (or service level) control mode. We consider SCC as a special operational mode that is activated across a pre-defined 1-hour period for a 24-hour simulation.

\begin{figure}[ht]
    \centering
    \includegraphics[width=0.725\textwidth]{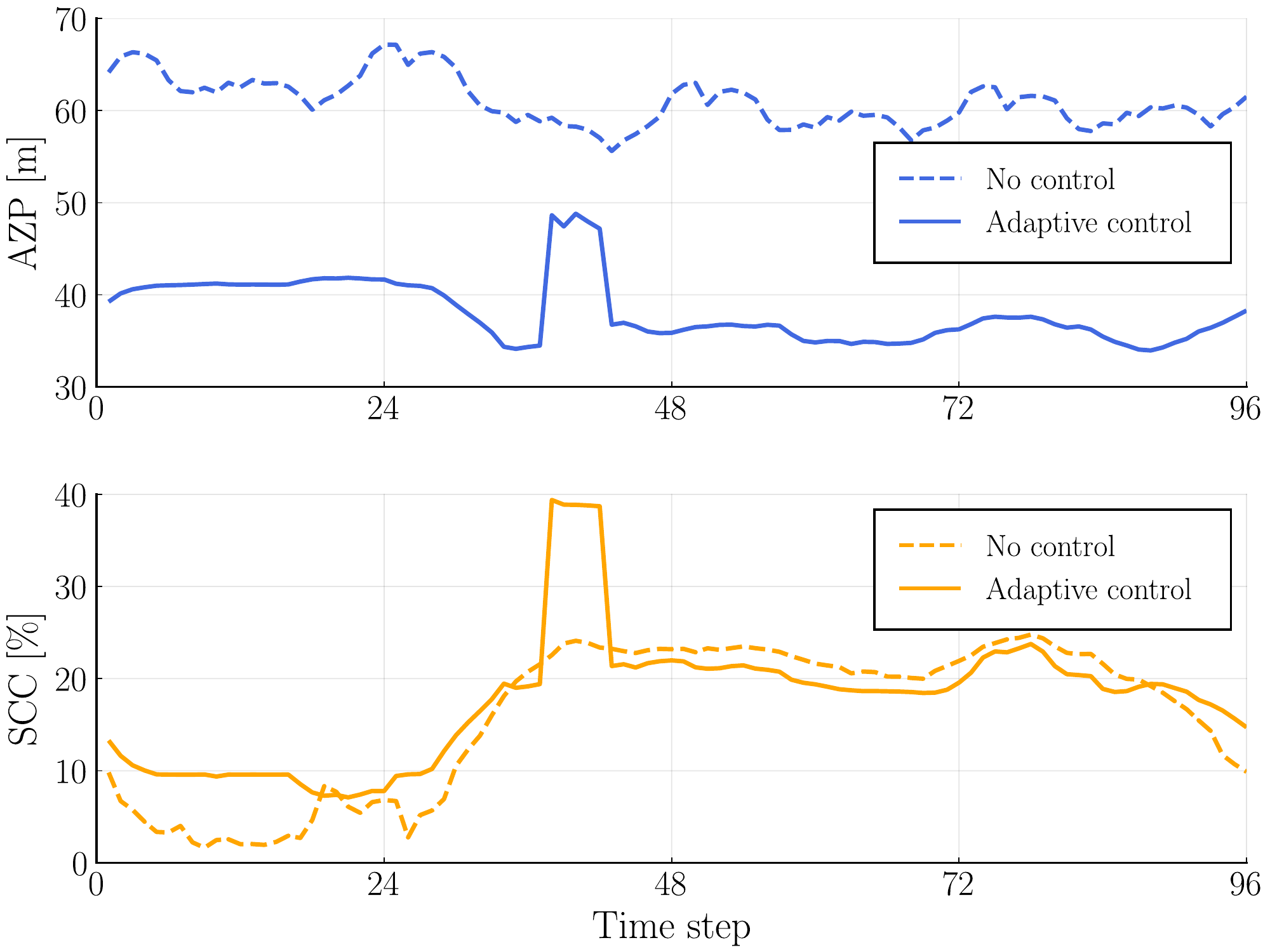}
    \caption{Time series of AZP-SCC adaptive control scheme during peak demand in \texttt{BWFLnet}.}
    \label{fig:adaptive_control}
\end{figure}

We demonstrate an example of preliminary results via a time series plot in \Cref{fig:adaptive_control} for \texttt{BWFLnet}. Here, SCC controls are optimized during the peak demand period (e.g. 09:30h to 10:30h), with AZP optimized across all other periods in the 24-hour simulation. As discussed in \Cref{sec:num_experiment_bi_dfc}, we optimize PCV actuators for both AZP and SCC control modes and AFV actuators for only the SCC mode. This adaptive control scheme enables the integration of pressure management controls, whilst facilitating a self-cleaning mode to reduce discolouration risk. Note that each time step in \Cref{fig:adaptive_control} represents a 15-minute period.

\begin{figure}[ht!]
    \centering
    \subfloat[\centering \label{fig:npv_cdf_modena}\texttt{Modena} ]{
    \includegraphics[width=0.44\textwidth]{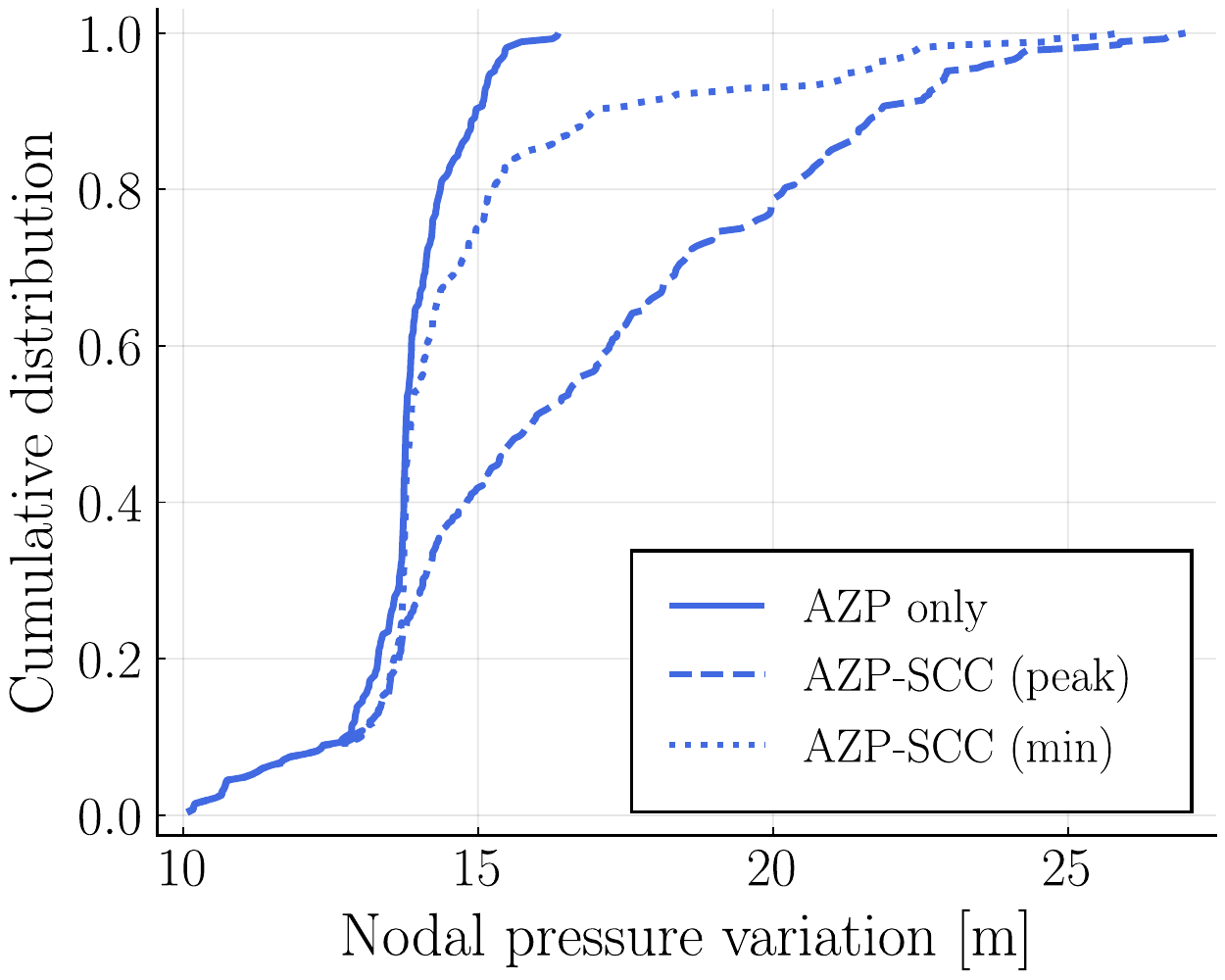}}
    \qquad
    \subfloat[\centering \label{fig:npv_cdf_bwfl}\texttt{BWFLnet} ]{
    \includegraphics[width=0.44\textwidth]{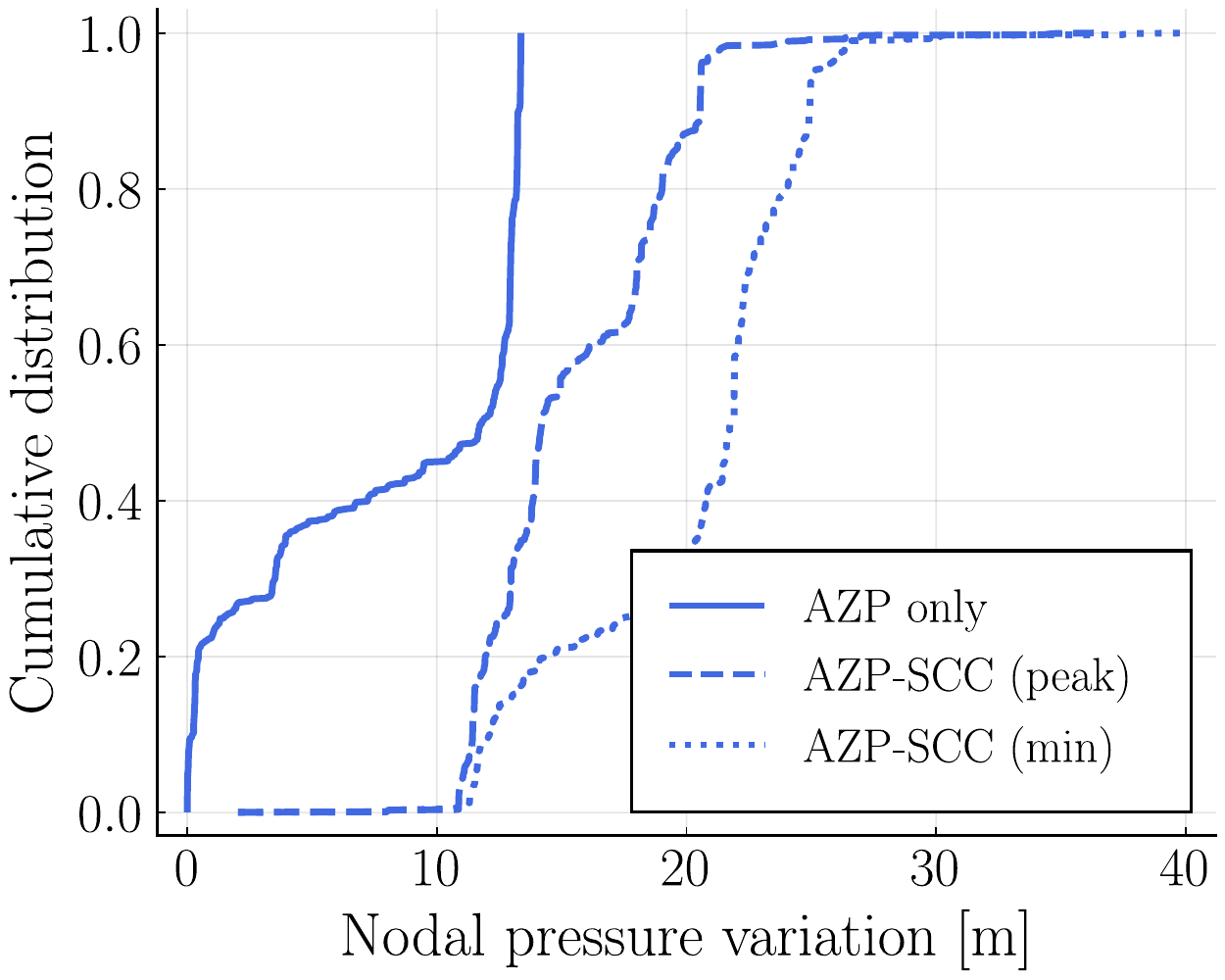}}
    \caption{Cumulative distribution plots of nodal pressure variation for different adaptive control scenarios.}
    \label{fig:npv_cdf}
\end{figure}

As expected, we observe a relatively large increase in AZP when activating optimal SCC controls. These conditions are commensurate with the SCC anchor point shown in \Cref{fig:bwfl_pareto}. The transition between control modes therefore result in a wider range in pressures across the network. While the instantaneous pressure changes would be smoothed by the settings of PID (Proportional-Integral-Derivative) controls at the device level, the network's physical infrastructure would still experience cyclic loading due to daily variations in pressure. Recalling the literature review in \Cref{sec:review_obj}, the frequency and range of nodal pressure variation has been established as a potential influence on pipe failure rates. Consequently, maximum pressure and pressure variability needs to be constrained within acceptable bounds when implementing such adaptive control schemes. To get a better understanding of the potential range in nodal pressure, we evaluate conditions resulting from different SCC activation periods. \Cref{fig:npv_cdf} compares three such control scenarios simulated in \texttt{Modena} and \texttt{BWFLnet}: (i) minimization of AZP only, (ii) adaptive AZP-SCC scheme with switch to SCC maximization during peak demand, and (iii) adaptive AZP-SCC scheme with switch to SCC maximization at the minimum demand period. These are conveyed through cumulative distribution plots of the range in nodal pressure variation. Note that the AZP-only control scenario represents a baseline condition for purposes of comparison.

For \texttt{Modena}, the activation of SCC controls at peak demand, relative to the AZP-only scenario, results in an increased pressure variation of 5~m to 10~m in approximately 60\% of model nodes. In comparison, we observe most nodes to have similar conditions between the AZP-only and minimum hour control scenarios, indicating less overall stress on the network. While not shown here for the sake of brevity, there is a modest decrease in SCC performance when implementing optimal SCC controls over the minimum demand hour compared to peak demand. This suggests an additional trade-off exists between SCC performance and nodal pressure variation. For \texttt{BWFLnet}, implementing optimal SCC controls at peak demand results in less pressure variation over the entire set of model nodes, compared to maximizing SCC over the minimum demand hour. Reviewing the optimal controls more closely, we note that flushing demands are higher at the minimum demand hour than in the peak demand scenario. Moreover, there are larger pressures at source nodes under minimum demand conditions, which results in a wider allowable pressure range in relation to the minimum regulatory pressure head. The difference in pressure variation for \texttt{BWFLnet} is therefore influenced by many network and operational specific factors. In any case, we indeed observe an increase in pressure variation compared to the AZP-only scenario, which is on the order of approximately 10 m to 15 m. We note that there is also a modest reduction in SCC performance in \texttt{BWFLnet} when optimal SCC controls are implemented over the minimum demand period. 

These findings demonstrate that increased nodal pressure variation arises from integrating AZP and SCC control modes. Since this represents a cyclic load on the network's physical infrastructure, future work might consider formulating an adaptive control problem constraining maximum pressure and pressure variability. However, this can be difficult to solve for large problem instances since the pressure variation metric adopted in \Cref{sec:review_obj_pressure} links hydraulic conditions between time steps. To this end, decomposition methods might be well suited to separate the large-scale problem into smaller sub-problems that can be solved efficiently. Furthermore, quantifying the impact of pressure variation on pipe burst rates is not well understood. Accordingly, further research is needed to consider a maximum allowable pressure variation in the implementation of such control problems.

\section{Conclusions}
\label{sec:conclusion}
We have investigated the problem of integrating optimal pressure management and self-cleaning operations in dynamically adaptive water distribution networks. The considered objectives are the minimization of average zone pressure (AZP) for reducing leakage and maximization of network self-cleaning capacity (SCC) for reducing discolouration risk. We first presented a review of existing single-objective AZP and SCC control strategies. This concluded with two current challenges concerning the integration of AZP and SCC objectives, which comprised the contributions of this paper.

We first formulated and solved the AZP-SCC bi-objective problem using the weighted sums scalarization method. By leveraging a previously developed convex heuristic, we approximated the Pareto fronts for two specific problems: (i) the design-for-control problem, where both locations and operational settings of PCV and AFV actuators were jointly optimized; and (ii) a control-only problem, where valve actuator settings were optimized for a fixed valve configuration corresponding to a hierarchical design strategy. The hierarchical design considered PCV actuators placed for the individual AZP problem and, following this solution, AFV actuators placed to augment SCC conditions. The resulting Pareto fronts suggest that significant improvements in SCC can be achieved for minimal trade-offs in AZP performance. Moreover, the hierarchical design approach yielded good quality solutions by varying actuator settings for the fixed valve locations. This is a valuable insight towards the integration of water quality for dynamically adaptive control in smart water networks. Future work should also consider minimizing the costs associated with valve installations.

We then demonstrated an adaptive control scheme to dynamically transition between AZP and SCC control modes. We formulated this control problem as a sequence of single-objective optimization problems, whereby AZP and SCC were considered the primary (master) and secondary (service level) control modes, respectively. Preliminary results highlighted the potential for adverse increases in nodal pressure variation arising from the transition between control modes. Although changes in network hydraulics are inevitable, future work might consider constraining pressure variation to limit the magnitude of the resulting cyclic load. Furthermore, the proposed adaptive control scheme could be implemented in a model predictive control (MPC) framework combined with higher resolution customer metering data. This could enable a shorter control horizon, reducing uncertainties associated with the self-cleaning period, which is assumed as a known 1-hour window in this work.

The presented design and control methods enhance the capabilities of water network operators to cost-effectively and optimally integrate pressure management and self-cleaning operations within dynamically adaptive water distribution networks.

\section*{Acknowledgments}
\noindent This work was supported by EPSRC (EP/P004229/1, Dynamically Adaptive and Resilient Water Supply Networks for a Sustainable Future), Bristol Water Plc, Analytical Technology (ATi), and Imperial College London’s Department of Civil and Environmental Engineering Skempton Scholarship.

\cleardoublepage

\bibliography{mybibfile}

\begin{thebibliography}{}

\bibitem[Abraham et~al., 2016]{ABRAHAM2016}
Abraham, E., Blokker, M., and Stoianov, I. (2016).
\newblock Network analysis, control valve placement and optimal control of flow
  velocity for self-cleaning water distribution systems.
\newblock {\em 18th Conference on Water Distribution System Analysis,
  WDSA2016}, pages 1--9.

\bibitem[Abraham et~al., 2018]{ABRAHAM2018}
Abraham, E., Blokker, M., and Stoianov, I. (2018).
\newblock Decreasing the discoloration risk of drinking water distribution
  systems through optimized topological changes and optimal flow velocity
  control.
\newblock {\em Journal of Water Resources Planning and Management},
  144(2):04017093.

\bibitem[Abraham et~al., 2019]{ABRAHAM2019}
Abraham, E., Pecci, F., and Stoianov, I. (2019).
\newblock A multi-objective framework for managing self-cleaning capacity and
  leakage: application to a real network model.
\newblock {\em 17th International Computing \& Control for the Water Industry
  Conference, 1-4 September 2019, Exeter, United Kingdom}, pages 1--2.

\bibitem[Abraham and Stoianov, 2015]{ABRAHAM2015}
Abraham, E. and Stoianov, I. (2015).
\newblock Sparse null space algorithms for hydraulic analysis of large-scale
  water supply networks.
\newblock {\em Journal of Hydraulic Engineering}, 142(3).

\bibitem[Araujo et~al., 2006]{ARAUJO2006}
Araujo, L., Ramos, H., and Coelho, S.~T. (2006).
\newblock Pressure control for leakage minimisation in water distribution
  systems management.
\newblock {\em Water Resources Management}, 20:133--149.

\bibitem[Armand et~al., 2017]{ARMAND2017}
Armand, H., Stoianov, I., and Graham, N. (2017).
\newblock A holistic assessment of discolouration processes in water
  distribution networks.
\newblock {\em Urban Water Journal}, 14(3):263--277.

\bibitem[Armand et~al., 2018]{ARMAND2018}
Armand, H., Stoianov, I., and Graham, N. (2018).
\newblock Impact of network sectorisation on water quality management.
\newblock {\em Journal of Hydroinformatics}, 20(2):424--439.

\bibitem[Belotti et~al., 2009]{BELOTTI2009}
Belotti, P., Lee, J., Liberti, L., Margot, F., and Wächter, A. (2009).
\newblock Branching and bounds tightening techniques for non-convex minlp.
\newblock {\em Optimization Methods and Software}, 24(4-5):597--634.

\bibitem[Bezanson et~al., 2017]{BEZANSON2017}
Bezanson, J., Edelman, A., Karpinski, S., and Shah, V.~B. (2017).
\newblock Julia: A fresh approach to numerical computing.
\newblock {\em SIAM Review}, 59(1):65--98.

\bibitem[Blokker et~al., 2010]{BLOKKER2010}
Blokker, E., Vreeburg, J., Schaap, P., and van Dijk, J. (2010).
\newblock The self-cleaning velocity in practice.
\newblock {\em Water Distribution System Analysis (WDSA)}, pages 187--199.

\bibitem[Blokker et~al., 2012]{BLOKKER2012}
Blokker, M., Vogelaar, H., Goos, K., and Vreeburg, J. (2012).
\newblock Using valve manipulation to manage discolouration risk in drinkiing
  water.
\newblock {\em Water Asset Management International}, 1:7--10.

\bibitem[Boxall and Saul, 2005]{BOXALL2005}
Boxall, J. and Saul, A. (2005).
\newblock Modeling discoloration in potable water distribution systems.
\newblock {\em Environmental Engineering}, 131(5):716--725.

\bibitem[Braga et~al., 2020]{BRAGA2020}
Braga, A.~S., Saulnier, R., Filion, Y., and Cushing, A. (2020).
\newblock Dynamics of material detachment from drinking water pipes under
  flushing conditions in a full-scale drinking water laboratory system.
\newblock {\em Urban Water Journal}, 17(8):745--753.

\bibitem[Bragalli et~al., 2012]{BRAGALLI2012}
Bragalli, C., Lodi, A., and D'Ambrosio, C. (2012).
\newblock On the optimal design of water distribution networks: a practical
  minlp approach.
\newblock {\em Optimization and Engineering}, 13(2):219--246.

\bibitem[Covelli et~al., 2016]{COVELLI2016}
Covelli, C., Cozzolino, L., Cimorelli, L., Morte, R.~D., and Pianese, D.
  (2016).
\newblock Optimal location and setting of prvs in wds for leakage minimization.
\newblock {\em Water Resources Management}, 30(5):1803--1817.

\bibitem[Creaco and Pezzinga, 2015]{CREACO2015}
Creaco, E. and Pezzinga, G. (2015).
\newblock Multiobjective optimization of pipe replacements and control valve
  installations for leakage attenuation in water distribution networks.
\newblock {\em Journal of Water Resources Planning and Management},
  141(3):04014059.

\bibitem[Dai and Li, 2014]{DAI2014}
Dai, P.~D. and Li, P. (2014).
\newblock Optimal localization of pressure reducing valves in water
  distribution systems by a reformulation approach.
\newblock {\em Water Resources Management}, 28(10):3057--3074.

\bibitem[D'Ambrosio et~al., 2015]{DAMBROSIO2015}
D'Ambrosio, C., Lodi, A., Wiese, S., and Bragalli, C. (2015).
\newblock Mathematical programming techniques in water network optimization.
\newblock {\em European Journal of Operational Research}, 243(3):774--788.

\bibitem[Das and Dennis, 1997]{DAS1997}
Das, I. and Dennis, J.~E. (1997).
\newblock A closer look at drawbacks of minimizing weighted sums of objectives
  for pareto set generation in multicriteria optimization problems.
\newblock {\em Structural Optimization}, 14:63--69.

\bibitem[{De Paola} et~al., 2017]{DEPAOLA2017}
{De Paola}, F., Galdiero, E., and Giugni, M. (2017).
\newblock Location and setting of valves in water distribution networks using a
  harmony search approach.
\newblock {\em Journal of Water Resources Planning and Management},
  143(6):04017015.

\bibitem[Dini and Asadi, 2020]{DINI2020}
Dini, M. and Asadi, A. (2020).
\newblock Optimal operational scheduling of available partially closed valves
  for pressure management in water distribution networks.
\newblock {\em Water Resources Management}, 34(8):2571--2583.

\bibitem[{Drinking Water Inspectorate}, 2020]{DWI2020}
{Drinking Water Inspectorate} (2020).
\newblock Drinking water safety guidance to health and water professionals.

\bibitem[Dunning et~al., 2017]{DUNNING2017}
Dunning, I., Huchette, J., and Lubin, M. (2017).
\newblock Jump: A modeling language for mathematical optimization.
\newblock {\em SIAM Review}, 59(2):295--320.

\bibitem[Eck and Mevissen, 2015]{ECK2015}
Eck, B. and Mevissen, M. (2015).
\newblock Quadratic approximations for pipe friction.
\newblock {\em Journal of Hydroinformatics}, 17(3):462--472.

\bibitem[Eck and Mevissen, 2013]{ECK2013}
Eck, B.~J. and Mevissen, M. (2013).
\newblock Fast non-linear optimization for design problems on water networks.
\newblock {\em World Environmental and Water Resources Congress 2013}, pages
  696--705.

\bibitem[Ehrgott, 2006]{EHRGOTT2006}
Ehrgott, M. (2006).
\newblock A discussion of scalarization techniques for multiple objective
  integer programming.
\newblock {\em Annals of Operations Research}, 147(1):343--360.

\bibitem[Ghaddar et~al., 2017]{GHADDAR2017}
Ghaddar, B., Claeys, M., Mevissen, M., and Eck, B.~J. (2017).
\newblock Polynomial optimization for water networks: Global solutions for the
  valve setting problem.
\newblock {\em European Journal of Operational Research}, 261(2):450--459.

\bibitem[{Gurobi Optimization}, 2022]{GUROBI2022}
{Gurobi Optimization} (2022).
\newblock Gurobi optimizer 9.5.0 reference manual.

\bibitem[HSL, 2021]{HSL2021}
HSL (2021).
\newblock A collection of fortran codes for large scale scientific computation.

\bibitem[Jara-Arriagada and Stoianov, 2021]{ARRIAGADA2021}
Jara-Arriagada, C. and Stoianov, I. (2021).
\newblock Pipe breaks and estimating the impact of pressure control in water
  supply networks.
\newblock {\em Reliability Engineering and System Safety}, 210:107525.

\bibitem[Jenks et~al., 2023]{JENKS2023}
Jenks, B., Pecci, F., and Stoianov, I. (2023).
\newblock Optimal design-for-control of self-cleaning water distribution
  networks using a convex multi-start algorithm.
\newblock {\em Water Research}, 231:119602.

\bibitem[Jowitt and Xu, 1989]{JOWITT1989}
Jowitt, P.~W. and Xu, C. (1989).
\newblock Optimal valve control in water-distribution networks.
\newblock {\em Journal of Water Resources and Planning Management},
  116(4):455--472.

\bibitem[Lambert and Fantozzi, 2010]{LAMBERT2010}
Lambert, A. and Fantozzi, M. (2010).
\newblock Recent developments in pressure management.

\bibitem[Larock et~al., 1999]{LAROCK1999}
Larock, B., Jeppson, R., and Watters, G. (1999).
\newblock {\em Hydraulics of Pipeline Systems}.
\newblock CRC press.

\bibitem[{Local Government Association and Water UK}, 2007]{FIRE2007}
{Local Government Association and Water UK} (2007).
\newblock National guidance document on the provision of water for fire
  fighting.

\bibitem[Machell and Boxall, 2014]{MACHELL2014}
Machell, J. and Boxall, J. (2014).
\newblock Modeling and field work to investigate the relationship between age
  and quality of tap water.
\newblock {\em Journal of Water Resources Planning and Management},
  140(9):04014020.

\bibitem[Mala-Jetmarova et~al., 2017]{MALA2017}
Mala-Jetmarova, H., Sultanova, N., and Savic, D. (2017).
\newblock Lost in optimisation of water distribution systems? a literature
  review of system operation.
\newblock {\em Environmental Modelling and Software}, 93:209--254.

\bibitem[Martínez-Codina et~al., 2016]{MARTINEZ2016}
Martínez-Codina, Castillo, M., González-Zeas, D., and Garrote, L. (2016).
\newblock Pressure as a predictor of occurrence of pipe breaks in water
  distribution networks.
\newblock {\em Urban Water Journal}, 13(7):676--686.

\bibitem[Messac et~al., 2003]{MESSAC2003}
Messac, A., Ismail-Yahaya, A., and Mattson, C.~A. (2003).
\newblock The normalized normal constraint method for generating the pareto
  frontier.
\newblock {\em Structural and Multidisciplinary Optimization}, 25(2):86--98.

\bibitem[Miettinen, 1998]{MIETTINEN1998}
Miettinen, K. (1998).
\newblock {\em Nonlinear Multiobjective Optimization}, volume~12.
\newblock Springer US.

\bibitem[Nerantzis et~al., 2020]{NERANTZIS2020}
Nerantzis, D., Pecci, F., and Stoianov, I. (2020).
\newblock Optimal control of water distribution networks without storage.
\newblock {\em European Journal of Operational Research}, 284(1):345--354.

\bibitem[Nerantzis and Stoianov, 2022]{NERANTZIS2022}
Nerantzis, D. and Stoianov, I. (2022).
\newblock Adaptive model predictive control for fire incidents in water
  distribution networks.
\newblock {\em Journal of Water Resources Planning and Management},
  148(2):04021102.

\bibitem[Nicolini and Zovatto, 2009]{NICOLINI2009}
Nicolini, M. and Zovatto, L. (2009).
\newblock Optimal location and control of pressure reducing valves in water
  networks.
\newblock {\em Journal of Water Resources and Planning Management},
  135(3):178--187.

\bibitem[Ofwat, 2022]{OFWAT2022}
Ofwat (2022).
\newblock Water company performance report 2021-2022.

\bibitem[Pecci et~al., 2017a]{PECCI2017a}
Pecci, F., Abraham, E., and Stoianov, I. (2017a).
\newblock Quadratic head loss approximations for optimisation problems in water
  supply networks.
\newblock {\em Journal of Hydroinformatics}, 19(4):493--506.

\bibitem[Pecci et~al., 2017b]{PECCI2017b}
Pecci, F., Abraham, E., and Stoianov, I. (2017b).
\newblock Scalable pareto set generation for multiobjective co-design problems
  in water distribution networks: a continuous relaxation approach.
\newblock {\em Structural and Multidisciplinary Optimization}, 55(3):857--869.

\bibitem[Pecci et~al., 2019]{PECCI2019}
Pecci, F., Abraham, E., and Stoianov, I. (2019).
\newblock Global optimality bounds for the placement of control valves in water
  supply networks.
\newblock {\em Optimization and Engineering}, 20(2):457--495.

\bibitem[Pecci et~al., 2021]{PECCI2021}
Pecci, F., Stoianov, I., and Ostfeld, A. (2021).
\newblock Relax-tighten-round algorithm for optimal placement and control of
  valves and chlorine boosters in water networks.
\newblock {\em European Journal of Operational Research}, 295:690--698.

\bibitem[Pecci et~al., 2022a]{PECCI2022a}
Pecci, F., Stoianov, I., and Ostfeld, A. (2022a).
\newblock Convex heuristics for optimal placement and operation of valves and
  chlorine boosters in water networks.
\newblock {\em Journal of Water Resources Planning and Management},
  148(2):04021098.

\bibitem[Pecci et~al., 2022b]{PECCI2022b}
Pecci, F., Stoianov, I., and Ostfeld, A. (2022b).
\newblock Optimal design-for-control of chlorine booster systems in water
  networks via convex optimization.
\newblock pages 1988--1993. IEEE.

\bibitem[Piller and van Zyl, 2014]{PILLER2014}
Piller, O. and van Zyl, J.~E. (2014).
\newblock Incorporating the favad leakage equation into water distribution
  system analysis.
\newblock {\em Procedia Engineering}, 89:613--617.

\bibitem[Reis et~al., 1997]{REIS1997}
Reis, L., Porto, R., and Chaudhry, F. (1997).
\newblock Optimal location of control valves in pipe networks by genetic
  algorithm.
\newblock {\em Journal of Water Resources Planning and Management},
  123(6):317--326.

\bibitem[Rezaei et~al., 2015]{REZAEI2015}
Rezaei, H., Ryan, B., and Stoianov, I. (2015).
\newblock Pipe failure analysis and impact of dynamic hydraulic conditions in
  water supply networks.
\newblock {\em Procedia Engineering}, 119(1):253--262.

\bibitem[Rossman et~al., 2020]{EPANET2.2}
Rossman, L., Woo, H., Tryby, M., Shang, F., Janke, R., and Haxton, T. (2020).
\newblock Epanet 2.2 user manual.

\bibitem[Ryan et~al., 2008]{RYAN2008}
Ryan, G., Mathes, P., Haylock, G., Jayaratne, A., Wu, J., Noui-Mehidi, N.,
  Grainger, C., and Nguyen, B. (2008).
\newblock Particles in water distribution systems.

\bibitem[Schwaller and van Zyl, 2015]{SCHWALLER2015}
Schwaller, J. and van Zyl, J.~E. (2015).
\newblock Modeling the pressure-leakage response of water distribution systems
  based on individual leak behavior.
\newblock {\em Journal of Hydraulic Engineering}, 141(5).

\bibitem[Todini and Pilati, 1988]{TODINI1988}
Todini, E. and Pilati, S. (1988).
\newblock A gradient algorithm for the analysis of pipe networks.
\newblock {\em Computer applications in water supply: vol. 1---systems analysis
  and simulation}, pages 1--20.

\bibitem[Ulanicki et~al., 2000]{ULANICKI2000}
Ulanicki, B., Bounds, P. L.~M., Rance, J.~P., and Reynolds, L. (2000).
\newblock Open and closed loop pressure control for leakage reduction.
\newblock {\em Urban Water}, 2:105--114.

\bibitem[Ulusoy et~al., 2022a]{ULUSOY2022b}
Ulusoy, A.~J., Mahmoud, H.~A., Pecci, F., Keedwell, E.~C., and Stoianov, I.
  (2022a).
\newblock Bi-objective design-for-control for improving the pressure management
  and resilience of water distribution networks.
\newblock {\em Water Research}, 222:118914.

\bibitem[Ulusoy et~al., 2020]{ULUSOY2020}
Ulusoy, A.~J., Pecci, F., and Stoianov, I. (2020).
\newblock An minlp-based approach for the design-for-control of resilient water
  supply systems.
\newblock {\em IEEE Systems Journal}, 14(3):4579--4590.

\bibitem[Ulusoy et~al., 2022b]{ULUSOY2022a}
Ulusoy, A.~J., Pecci, F., and Stoianov, I. (2022b).
\newblock Bi-objective design-for-control of water distribution networks with
  global bounds.
\newblock {\em Optimization and Engineering}, 23(1):527--577.

\bibitem[Vairavamoorthy and Lumbers, 1998]{VAIRA1998}
Vairavamoorthy, K. and Lumbers, J. (1998).
\newblock Leakage reduction in water distribution systems: Optimal valve
  control.
\newblock {\em Journal of Hydraulic Engineering}, 124(11):1146--1154.

\bibitem[{van Summeren} and Blokker, 2017]{VANSUMMEREN2017}
{van Summeren}, J. and Blokker, M. (2017).
\newblock Modeling particle transport and discoloration risk in drinking water
  distribution networks.
\newblock {\em Drinking Water Engineering and Science}, 10(2):99--107.

\bibitem[Vreeburg et~al., 2009]{VREEBURG2009}
Vreeburg, J.~H., Blokker, E.~J., Horst, P., and van Dijk, J.~C. (2009).
\newblock Velocity-based self-cleaning residential drinking water distribution
  systems.
\newblock {\em Water Science and Technology: Water Supply}, 9(6):635--641.

\bibitem[Waldron et~al., 2020]{WALDRON2020}
Waldron, A., Pecci, F., and Stoianov, I. (2020).
\newblock Regularization of an inverse problem for parameter estimation in
  water distribution networks.
\newblock {\em Journal of Water Resources and Planning Management},
  146(9):04020076.

\bibitem[Wright et~al., 2015]{WRIGHT2015}
Wright, R., Abraham, E., Parpas, P., and Stoianov, I. (2015).
\newblock Control of water distribution networks with dynamic dma topology
  using strictly feasible sequential convex programming.
\newblock {\em Water Resources Research}, 51(12):9925--9941.

\bibitem[Wright et~al., 2014]{WRIGHT2014}
Wright, R., Stoianov, I., Parpas, P., Henderson, K., and King, J. (2014).
\newblock Adaptive water distribution networks with dynamically reconfigurable
  topology.
\newblock {\em Journal of Hydroinformatics}, 16(6):1280--1301.

\bibitem[Wächter and Biegler, 2006]{WACHTER2006}
Wächter, A. and Biegler, L.~T. (2006).
\newblock On the implementation of an interior-point filter line-search
  algorithm for large-scale nonlinear programming.
\newblock {\em Mathematical Programming}, 106(1):25--57.

\end{thebibliography}

\end{document}